\pdfoutput=1
\documentclass[11pt]{article}

\usepackage[T1]{fontenc}
\usepackage[utf8]{inputenc}
\usepackage{lmodern}
\usepackage[margin=1in]{geometry}
\usepackage{amsmath,amssymb,amsthm,mathtools}
\usepackage{physics}
\usepackage{bm}
\usepackage{microtype}
\usepackage{enumitem}
\usepackage{graphicx}
\usepackage{float}
\usepackage{tikz}
\usetikzlibrary{arrows.meta,positioning}
\usepackage{algorithm}
\usepackage{algorithmic}
\usepackage{booktabs}
\usepackage{subcaption}
\usepackage{siunitx}
\usepackage{tabularx}
\usepackage{makecell}
\usepackage[numbers,sort&compress]{natbib}
\usepackage{hyperref}
\hypersetup{
  colorlinks=true,
  linkcolor=blue,
  citecolor=blue,
  urlcolor=blue,
  pdftitle={Quaternionic and Octonionic Frameworks for Quantum Computation: Mathematical Structures, Models, and Fundamental Limitations},
  pdfauthor={Johan Heriberto Rua Munoz; Jorge Eduardo Mahecha Gomez; Santiago Pineda Montoya}
}

\newcommand{\HH}{\mathbb{H}}
\newcommand{\OO}{\mathbb{O}}
\newcommand{\C}{\mathbb{C}}

\theoremstyle{plain}
\newtheorem{Theorem}{Theorem}
\newtheorem{Lemma}{Lemma}
\newtheorem{Corollary}{Corollary}
\newtheorem{Proposition}{Proposition}
\theoremstyle{remark}
\newtheorem{Remark}{Remark}

\title{Quaternionic and Octonionic Frameworks for Quantum Computation: Mathematical Structures, Models, and Fundamental Limitations}

\author{
Johan Heriberto R\'{u}a Mu\~{n}oz$^{1}$,
Jorge Eduardo Mahecha G\'{o}mez$^{1}$, and
Santiago Pineda Montoya$^{2,*}$\\[0.6em]
\small $^{1}$ Physics Institute, Exact and Natural Sciences Faculty, Antioquia University, Medell\'{i}n 050010, Colombia\\
\small $^{2}$ Mathematics Institute, Sciences Faculty, National University of Colombia, Medell\'{i}n 050034, Colombia\\
\small $^{*}$ Corresponding author: \href{mailto:sapinedamo@unal.edu.co}{sapinedamo@unal.edu.co}\\
\small Emails: \href{mailto:heriberto.rua@udea.edu.co}{heriberto.rua@udea.edu.co};
\href{mailto:jorge.mahecha@udea.edu.co}{jorge.mahecha@udea.edu.co}
}

\date{\small Published in \textit{Quantum Reports}, 7(4), 55 (2025).\\
DOI: \href{https://doi.org/10.3390/quantum7040055}{10.3390/quantum7040055}
\quad Published: 26 November 2025}

\begin{document}
\maketitle

\begin{abstract}
We develop detailed quaternionic and octonionic frameworks for quantum computation grounded on normed division algebras. Our central result is to prove the polynomial computational equivalence of quaternionic and complex quantum models: Computation over $\mathbb{H}$ is polynomially equivalent to the standard complex quantum circuit model and hence captures the same complexity class BQP up to polynomial reductions. Over $\mathbb{H}$, we construct a complete model---quaternionic qubits on right $\mathbb{H}$-modules with quaternion-valued inner products, unitary dynamics, associative tensor products, and universal gate sets---and establish polynomial equivalence with the standard complex model; routes for implementation at fidelities exceeding $99\%$ via pulse-level synthesis on current hardware are discussed. Over $\mathbb{O}$, non-associativity yields path-dependent evolution, ambiguous adjoints/inner products, non-associative tensor products, and possible failure of energy conservation outside associative sectors. We formalize these obstructions and systematize four mitigation strategies: Confinement to associative subalgebras, $G_{2}$-invariant codes, dynamical decoupling of associator terms, and a seven-factor algebraic decomposition for gate synthesis. The results delineate the feasible quaternionic regime from the constrained octonionic landscape and point to applications in symmetry-protected architectures, algebra-aware simulation, and hypercomplex learning.
\end{abstract}

\noindent\textbf{Keywords:} quantum computation; quaternions; octonions; normed division algebras; non-associative quantum mechanics; universal gates; $G_2$ symmetry; dynamical decoupling.

\vspace{1em}

\section{Introduction}
\textls[-20]{Quantum computation is typically formulated on complex Hilbert spaces. The existence of exactly four real normed division algebras, $\mathbb{R},\mathbb{C},\mathbb{H},\mathbb{O}$, motivates systematic extensions to $\mathbb{H}$ (associative, non-commutative) and $\mathbb{O}$ (alternative, non-associative) \cite{Hurwitz1898,Baez2002}. We present a comprehensive treatment of both directions, revealing a sharp contrast: A fully consistent quaternionic model with computational equivalence to the complex standard \cite{Finkelstein1962,Adler1995,RuaThesis2025}, and structural obstructions in the octonionic case \cite{Baez2002,RuaThesis2025}.}

Main contributions.

We explicitly distinguish between background material (review and synthesis of known results) and contributions that, to the best of our knowledge, are original to \linebreak  this work.
\begin{itemize}[leftmargin=*,nosep]
  \item { Review and unification of quaternionic frameworks.
} We collect and streamline existing results on quaternionic quantum mechanics and computation, including the circuit model of Fernández and Schneeberger~\cite{Fernandez2008}, the simulation of quaternionic dynamics and channels on complex Hilbert spaces~\cite{Graydon2013,Gantner2017}, and classical treatments of quaternionic quantum theory~\cite{Finkelstein1962,Adler1995}. We recast these in a right-$\HH$-module language aligned with the modern quantum-information toolkit (states, channels, POVMs, and tensor products).
  \item {Original: encoded $\mathbb{H}$-based model of quantum computation.} Building on this background, we formulate a complete $\HH$-based circuit model with well-defined quaternionic qubit registers, gates, measurements, and tensor products, together with an explicit embedding into standard complex qubit circuits via the map $\iota\!:\HH\hookrightarrow M_2(\C)$. Our construction yields a polynomial-time, ancilla-efficient simulation of any quaternionic process by complex circuits, extending the circuit-level result of~\cite{Fernandez2008} to general channels and measurement patterns and making the computational equivalence with BQP fully explicit (see Theorem~\ref{thm:eq} and Proposition~\ref{prop:complexH}).
  \item {Review of octonionic proposals.} We summarize relevant prior work on octonions in physics and quantum information, in particular the measurement-only scheme for universal quantum computation based on octonionic equiangular projections~\cite{Freedman2018} and earlier studies of octonionic structures in quantum theory~\cite{Baez2002}.
  \item {Original: systematic obstruction and mitigation analysis for $\mathbb{O}$-based models.} We show that promoting $\OO$ to a gate-based computational model faces three structural obstructions (loss of associative tensor products, path-dependent dynamics, and nonlocal measurement postulates), and we organize possible workarounds into four concrete strategies (M1)--(M4): associative confinement, $G_2$-invariant coding, associator-suppressing dynamical decoupling, and seven-factor synthesis. We analyse their resource overheads and limits in terms of qubit counts, control complexity, and fault-tolerance, clarifying in what sense octonionic frameworks can be realised only as encoded or effective complex-qubit models.
\end{itemize}
\vspace{3pt}

Context and motivation.

Hurwitz’s theorem isolates $\mathbb{R},\mathbb{C},\mathbb{H}$, and $\mathbb{O}$ as the only finite-dimensional real normed division algebras, providing a canonical ladder of algebraic structure via the Cayley--Dickson construction. Within this ladder, $\mathbb{H}$ preserves associativity (though not commutativity), whereas $\mathbb{O}$ forfeits associativity but remains alternative \cite{Hurwitz1898,Baez2002}. These properties are not merely algebraic curiosities: Associativity underpins the conventional tensorial wiring of multipartite quantum systems, while alternativity constrains but does not fully restore familiar quantum identities. Against this backdrop, we develop a side-by-side account that makes precise which elements of the standard quantum formalism survive over $\mathbb{H}$ and which fail over $\mathbb{O}$ \cite{Adler1995,RuaThesis2025}.

Quaternionic model in brief.

Over $\mathbb{H}$, we adopt right $\mathbb{H}$-module “Hilbert” spaces with quaternion-valued inner products that are anti-linear in the first slot and linear in the second; observables are self-adjoint, and dynamics is unitary in the usual sense \cite{Adler1995}. Crucially, tensor products are associative, so composition of subsystems is well-defined and independent of parenthesization. Via the multiplicative embedding $\iota:\mathbb{H}\hookrightarrow M_2(\mathbb{C})$, where $\iota(a+bi+cj+dk)=\begin{pmatrix} a+bi & c+di\\ -c+di & a-bi \end{pmatrix}$, unit quaternions act as $SU(2)$ rotations, allowing standard universal gate sets (single-qubit rotations plus a fixed entangler) to lift to the quaternionic setting. Using this embedding, we establish polynomial-time simulations in both directions, showing computational equivalence between the quaternionic and complex models and enabling pulse-level realizations aligned with native $SU(2)$ controls on present-day platforms \cite{Finkelstein1962,Adler1995,RuaThesis2025}.

Octonionic landscape in brief.

In stark contrast, non-associativity in $\mathbb{O}$ produces path-dependent evolution (propagators acquire associator corrections), makes adjoints and inner products order-sensitive, obstructs a canonical multipartite tensor structure, and can undermine standard conservation arguments outside associative sectors \cite{Baez2002,RuaThesis2025}. We formalize these obstructions and delineate their impact on circuit composition and physical interpretability, identifying the precise points at which alternativity is insufficient to recover standard quantum behavior.

Mitigation palette and scope.

While a fully “native’’ octonionic architecture appears structurally constrained, we systematize four principled mitigation routes: (M1) confinement to embedded quaternionic subalgebras, (M2) $G_2$-invariant coding that detects associator syndromes, (M3) dynamical decoupling tailored to cancel coherent associator contributions at low orders, and (M4) a seven-factor algebraic synthesis in associative sectors. Each comes with quantifiable overheads (depth, calibration, or restricted reachable sets) and is best understood as a means to simulate or protect limited octonionic behavior rather than to scale it without constraints \cite{RuaThesis2025}.

Paper organization.

Section~\ref{sec:methods} reviews the complex formalism and algebraic preliminaries. Section~\ref{sec:methods} develops the $\mathbb{H}$ and $\mathbb{O}$ frameworks. Section~\ref{sec:results} presents the core theorems: Associative tensor structure and universality over $\mathbb{H}$; path dependence, energy-balance issues, and the absence of a canonical tensor product over $\mathbb{O}$; together with mitigation strategies and their trade-offs. The discussion and conclusions synthesize implications for complexity, compilation, and near-term demonstrations \cite{Adler1995,Baez2002,RuaThesis2025}.


\subsection{Extended Related Work}
Background on quantum computation, topological phases, and materials that motivate robust encodings (see e.g., \cite{Fazekas1996,Dray2000,Ashtekar1994,Altland1997,Bott1970,Born1928,Gullion1990,Haeberlen1968,Viola1999,Guo2024,DeLeo2012,BayroCorrochano2018,Alicata2018,Ando2015,Asboth2016}).

Mathematical foundations: normed division algebras, exceptional symmetries, and representation-theoretic tools (see e.g., \cite{Bernevig2006,Bernevig2013,Bistritzer2010,Bistritzer2011,Brihuega2012,Brzezicki2019,Bychkov1984,Castro2007,Cayssol2013,Chaikin1995,Chang2013,Chang1995,Chen2019,Chen2020,Chittari2019}).

Control and error mitigation, with emphasis on average-Hamiltonian theory and dynamical decoupling (see e.g., \cite{Chiu2016,Chiu2014,Choi2019,DasSarma2011,Davila2014,deHeer2010,deLeseleuc2019,Deutsch1992,Duan2014,Dyakonov1971a,Dyakonov1971b,Eguchi1980,Epple1998,Ezawa2018,Ezawa2019}).

Complexity, universality, and emergent applications across condensed matter and quantum information (see e.g., \cite{Fang2003,Franz2013,Fu2006,Fu2007,Fu2020,Geim2013,Ghatak2019,Gong2014,Goss2022,Guo2009,Haldane1988,Hall2015,Hall1879,Hall1881,Halperin1982}).

\subsection{Introduction: Bridging Theory and Experiment—Additional Details}

The transition from an abstract mathematical formalism to a concrete physical implementation is the most challenging and critical step in the development of any new computational paradigm. This section serves as a bridge between the theoretical work and the practical considerations of experimental realization. We will first frame the core experimental challenges posed by the quaternionic and octonionic models and then provide a brief outline of the topics to be covered in this work \cite{RuaThesis2025}.

\subsubsection*{From
 Abstract Algebras to Physical Systems}

The experimental questions for realizing quaternionic and octonionic quantum computation are fundamentally different, reflecting the deep structural differences between the two algebras.

For quaternions, the path to implementation is relatively straightforward. As we proved in Chapter 3
 \cite{RuaThesis2025}, the group of logical operations on a single quaternionic qubit is isomorphic to $SU(2)$, the same group that governs standard complex qubits. This crucial fact implies that no fundamentally new hardware or physics is required to perform quaternionic quantum computation. A quaternionic gate is simply a specific type of rotation on the Bloch hypersphere, which can be decomposed and implemented as a sequence of control pulses on any universal quantum computer. The challenge is therefore not one of discovery, but of control: ``How can we efficiently engineer the laser or microwave pulses required to execute a desired quaternionic transformation on a standard qubit or qudit system?'' \cite{Adler1995}.

For octonions, the situation is far more speculative and challenging. A ``true'' octonionic quantum computer, where the elementary physical components themselves obeyed octonionic algebra, would imply the discovery of new physics beyond the Standard Model. A more practical, and the only currently conceivable, path is through simulation or embedding. As established in {Chapter 4} \cite{RuaThesis2025}, a single logical octo-qubit must be encoded into a larger register of standard qubits (a minimum of three) or the higher energy levels \linebreak  of a qudit.

The experimental challenge for octonions is therefore twofold:
\begin{enumerate}
    \item \textbf{The Challenge of Scale and Coherence:
} Can we build, initialize, and maintain the coherence of the complex, multi-qubit entangled states required to represent even a single octo-qubit?
    \item \textbf{The Challenge of Control and Fidelity:} Can we implement the sophisticated, multi-gate sequences required for the mitigation strategies (e.g., $G_2$-symmetric codes or dynamical decoupling) with a fidelity high enough to overcome the intrinsic \linebreak  associator errors?
\end{enumerate}

In essence,
 for quaternions, we ask how to express a new logic on existing hardware. For octonions, we ask if we can build sufficiently complex and precise hardware to tame a fundamentally different and more challenging form of logic \cite{Freedman2018}.

\section{Materials and Methods}\label{sec:methods}

\subsection{Standard Quantum Computation over \texorpdfstring{$\mathbb{C}$}{C}}
We briefly recall the standard (complex) formalism to fix notation. A finite-dimensional quantum system is modeled by a complex Hilbert space $\mathcal{H}\cong \mathbb{C}^{d}$ with an inner product $\langle \cdot,\cdot\rangle$ linear in the second argument and anti-linear in the first. \emph{States} are density operators $\rho\in\mathcal{D}(\mathcal{H})$ (positive semidefinite, trace one); pure states are unit vectors $\ket{\psi}\in \mathcal{H}$ modulo global phase. \emph{Closed} system dynamics is generated by a (time-dependent) Hamiltonian $H(t)=H(t)^\dagger$ through the Schrödinger equation:
\[
i\hbar \,\frac{d}{dt}\ket{\psi(t)}=H(t)\ket{\psi(t)},\qquad \ket{\psi(t)}=U(t)\ket{\psi(0)},\ U(t)\in U(d).
\]
Open
 evolutions are completely positive, trace-preserving (CPTP) maps $\Phi(\rho)=\sum_k E_k \rho E_k^\dagger$ (Kraus form) with $\sum_k E_k^\dagger E_k=\mathbb{I}$. Measurements are described by POVMs $\{M_j\}_{j}$ with $M_j\ge 0$, $\sum_j M_j=\mathbb{I}$; outcome $j$ occurs with probability $p_j=\Tr(M_j\rho)$ and post-measurement state $\rho_j=K_j\rho K_j^\dagger/p_j$ for an instrument $\{K_j\}_j$ s.t.\ $K_j^\dagger K_j=M_j$.

For one qubit, the Bloch representation writes $\rho=\tfrac12(\mathbb{I}+\vec r\cdot\vec \sigma)$, $\|\vec r\|\le 1$, with Pauli matrices $\sigma_x,\sigma_y,\sigma_z$ and rotations $R_{\hat n}(\theta)=\exp(-i\theta\, \hat n\cdot \vec \sigma/2)\in SU(2)$. Universal gate sets include any dense set of single-qubit rotations plus an entangling two-qubit gate \mbox{(e.g., CNOT)}. For reference,
\[
\mathrm{CNOT}=
\begin{pmatrix}
1&0&0&0\\
0&1&0&0\\
0&0&0&1\\
0&0&1&0
\end{pmatrix},
\]
acting on the computational basis $\{\ket{00},\ket{01},\ket{10},\ket{11}\}$. Multipartite systems are formed via the tensor product $\mathcal{H}_{AB}=\mathcal{H}_A\otimes \mathcal{H}_B$, with separable and entangled states defined as usual. These elements will serve as the comparison baseline for the quaternionic and octonionic extensions developed below \cite{Adler1995,RuaThesis2025}.

\subsection{Normed Division Algebras}\label{subsec:NDAs}
A (real) \emph{normed division algebra} $(A,\|\cdot\|)$ is a finite-dimensional real algebra with unit $1$, equipped with a norm such that
\[
\|xy\|=\|x\|\,\|y\|\qquad \forall\,x,y\in A,
\]
and where every nonzero element is invertible. A classical result due to Hurwitz shows that there are exactly four such algebras:
\[
\mathbb{R},\qquad \mathbb{C},\qquad \mathbb{H},\qquad \mathbb{O},
\]
with $\dim_{\mathbb{R}}=1,2,4,8$, respectively \cite{Hurwitz1898,Baez2002}. They arise from the Cayley--Dickson construction:
\[
\mathbb{R}\ \xrightarrow{\text{double}}\ \mathbb{C}\ \xrightarrow{\text{double}}\ \mathbb{H}\ \xrightarrow{\text{double}}\ \mathbb{O},
\]
and carry an involutive conjugation $x\mapsto \bar x$, a quadratic norm $N(x)=x\bar x=\bar x x=\|x\|^2\in \mathbb{R}_{\ge 0}$, and a bilinear inner product $\langle x,y\rangle := \tfrac12\big(N(x+y)-N(x)-N(y)\big)$.

Algebraically, $\mathbb{R}$ and $\mathbb{C}$ are commutative and associative; $\mathbb{H}$ is non-commutative but associative; $\mathbb{O}$ is \emph{non-associative} but \emph{alternative}, i.e., the subalgebra generated by any two elements is associative. The \emph{associator}
\[
[x,y,z]:=(xy)z - x(yz)
\]
vanishes identically in $\mathbb{R},\mathbb{C},\mathbb{H}$ and controls the failure of associativity in $\mathbb{O}$. Octonions obey the Moufang identities (e.g., $(xy)(zx)=x(yz)x$), which are often used to tame non-associativity in structured calculations \cite{Baez2002}.

Cayley--Dickson product.

Writing $\mathbb{O}=\mathbb{H}\oplus \mathbb{H}e_6$ as a real vector space means that every octonion $x\in\mathbb{O}$ can be written uniquely as follows:
\[
x = a + b e_6,\qquad a,b\in\mathbb{H},
\]
so that we may freely identify $x$ with the pair $(a,b)\in\mathbb{H}\oplus\mathbb{H}e_6$. With this convention, the Cayley--Dickson product for $(a,b),(c,d)\in\mathbb{H}\oplus\mathbb{H}e_6$ reads
\[
(a,b)\cdot(c,d)=(ac-\bar{d}\,b,\ dc+b\,\bar c),
\]
which makes explicit both the embedding of quaternionic subalgebras inside $\mathbb{O}$ and how non-associativity emerges beyond $\mathbb{H}$. This notation will be used repeatedly when discussing associative triples and confinement strategies.

\subsection{Quaternionic and Octonionic Preliminaries}\label{subsec:prelimHO}
The hypercomplex extensions considered here replace the complex scalar field by either $\mathbb{H}$ or $\mathbb{O}$, together with compatible notions of inner product, adjoint, dynamics,  and~composition.

Notation on $SU(2)$ and context.

Unless explicitly stated otherwise, $SU(2)$ denotes the standard group of $2\times 2$ \emph{complex} special unitaries acting on a single complex qubit. When we speak of ``quaternionic gates realised as $SU(2)$ operations'' or of ``native $SU(2)$ control'' on hardware, we always mean the image of quaternionic unitaries under the embedding $\iota:\HH\hookrightarrow M_2(\C)$, acting on an encoded complex qubit.

Quaternionic Hilbert modules.

A (right) quaternionic Hilbert module $(\mathcal{H}_{\mathbb{H}},\langle \cdot,\cdot\rangle)$ is a right $\mathbb{H}$-module endowed with an $\mathbb{H}$-valued inner product satisfying the following:
\begin{enumerate}
\item $\langle \psi,\phi a+\chi b\rangle=\langle \psi,\phi\rangle a+\langle \psi,\chi\rangle b$ \quad (right linearity);
\item $\langle \phi,\psi\rangle=\overline{\langle \psi,\phi\rangle}$ \quad (quaternionic Hermiticity);
\item $\langle \psi,\psi\rangle \in \mathbb{R}_{\ge 0}$ and vanishes iff $\psi=0$ \quad (positivity).
\end{enumerate}

{States} (\emph{quaternionic qubits}) have the form $\ket{\phi}=q_1\ket{0}+q_2\ket{1}$ with $q_i\in \mathbb{H}$ and\linebreak   $\|q_1\|^2+\|q_2\|^2=1$. Born probabilities are real: $P(\psi|\phi)=\|\langle \psi,\phi\rangle\|^2\in\mathbb{R}_{\ge 0}$. Linear operators are right $\mathbb{H}$-linear maps $A:\mathcal{H}_{\mathbb{H}}\to \mathcal{H}_{\mathbb{H}}$ with an adjoint $A^\dagger$ determined by $\langle \psi, A\phi\rangle=\langle A^\dagger\psi,\phi\rangle$. Unitary dynamics $U^\dagger U=\mathbb{I}$ is well-defined, and---crucially—tensor products over $\mathbb{H}$   are \emph{associative}:
\[
(\mathcal{H}_1\otimes \mathcal{H}_2)\otimes \mathcal{H}_3 \cong \mathcal{H}_1\otimes(\mathcal{H}_2\otimes \mathcal{H}_3),
\]
with inner product on simple tensors $\langle u_1\!\otimes\!u_2, v_1\!\otimes\!v_2\rangle=\langle u_1,v_1\rangle\,\langle u_2,v_2\rangle$. This guarantees a canonical multipartite structure and permits the usual notions of separability and entanglement (Theorem~\ref{thm:assocH} in the main text).

Matrix embedding \texorpdfstring{$\mathbb{H}\hookrightarrow M_2(\mathbb{C})$}{H->M2(C)}.

The multiplicative embedding $\iota:\ \mathbb{H}\to M_2(\mathbb{C})$, sends unit quaternions to $SU(2)$ and preserves norms up to a factor ($\|\iota(q)\|_F^2=2\|q\|^2$). It allows one to \emph{simulate} quaternionic circuits in the complex model with polynomial overhead (and conversely embed complex circuits into the quaternionic setting), and implies that native $SU(2)$ control translates directly into single-quaternionic qubit gate synthesis (see Appendix~\ref{app:HtoC}).

Octonionic modules and the associator.

Attempting to transplant the above structure to $\mathbb{O}$ runs into foundational obstacles controlled by the associator $[x,y,z]$:
\begin{itemize}
\item \emph{Adjoint/inner product ambiguity:} Without associativity, the defining relation $\langle \psi, A\phi\rangle=\langle A^\dagger\psi,\phi\rangle$ becomes parenthesization-dependent unless dynamics and observables are confined to associative subalgebras.
\item \emph{Non-associative tensor products:} For right $\mathbb{O}$-modules $\mathcal{H}_A,\mathcal{H}_B,\mathcal{H}_C$, the identification $(\mathcal{H}_A\otimes \mathcal{H}_B)\otimes \mathcal{H}_C \cong \mathcal{H}_A\otimes(\mathcal{H}_B\otimes \mathcal{H}_C)$ is not canonical in general; different parenthesizations yield inequivalent scalar actions. Thus, there is no \emph{canonical} multipartite structure unless one restricts to an associative (quaternionic) sector.
\item \textls[-15]{\emph{Path dependence in evolution:} For $H(t)\in \mathbb{O}$, the propagator $U(t)=\mathcal{P}\exp[-\tfrac{i}{\hbar}\int_0^t H(\tau)d\tau]$ acquires higher-order terms involving associators in Magnus-like expansions.} Reordering control slices can change $U(t)$ even for the same time-averaged $H$, leading to intrinsic path dependence outside associative sectors.
\end{itemize}

These features motivate the mitigation palette developed later (associative confinement, $G_2$-invariant coding, associator-aware decoupling, and seven-factor synthesis).

\subsection{Structure of Octonions and the Fano Plane}\label{subsec:fano}
\begin{figure}[H] \vspace{-6pt}
  \hspace{-28pt}\includegraphics[width=0.58\linewidth]{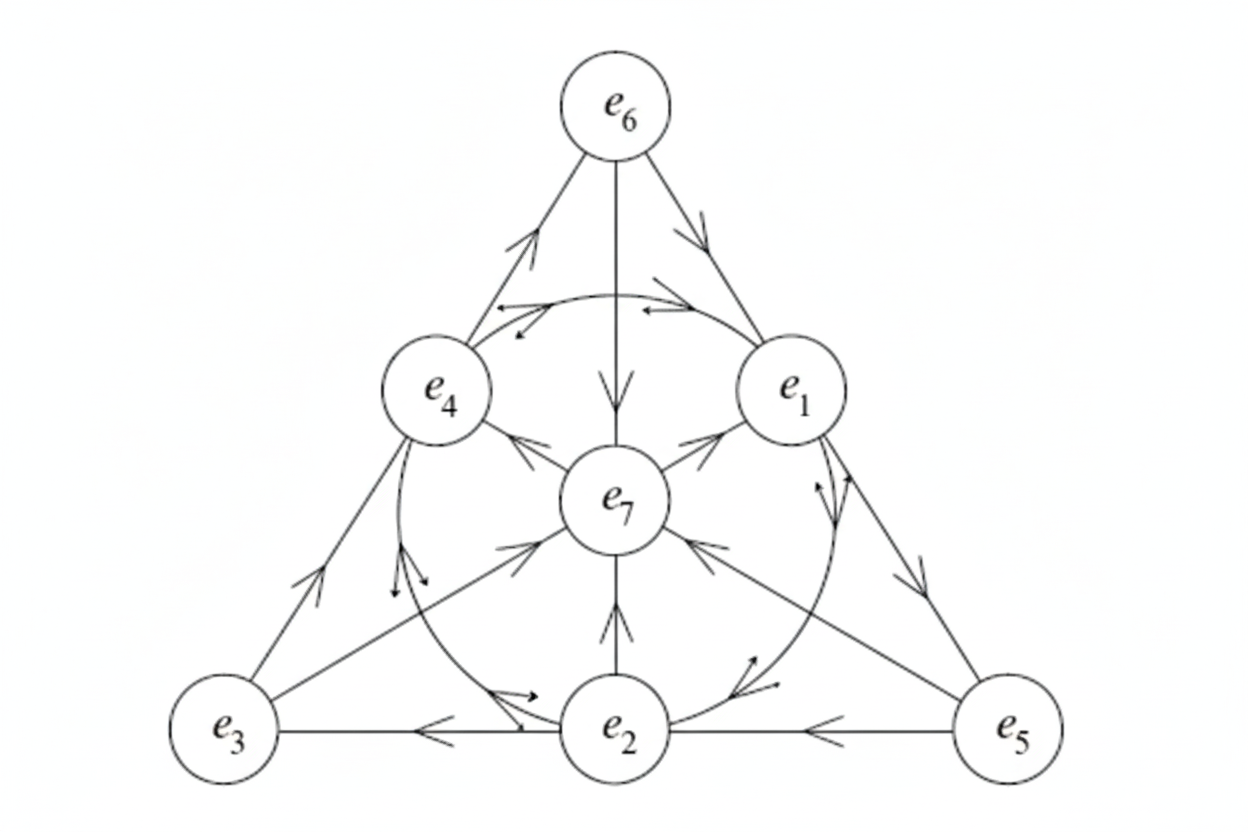}
  \caption{Fano
 plane representation of octonion multiplication triples. Each oriented line encodes an associative copy of $\mathbb{H}$ inside $\mathbb{O}$.}
  \label{fig:fano}
\end{figure}

Let $\{1,e_1,\dots,e_7\}$ be the standard octonion basis with $e_k^2=-1$ and $e_i e_j = \pm e_k$ according to the oriented triples encoded by the Fano plane (See Figure \ref{fig:fano}). Any three units on a directed line multiply associatively and obey $e_i e_j=e_k$, $e_j e_k=e_i$, $e_k e_i=e_j$ (and the negatives for reversed order). This diagrammatic calculus is indispensable for tracking when associators vanish and for designing sequences that remain within (or near) \mbox{associative triples.}

\subsection{Geometric and Group-Theoretic Side Notes}
The unit spheres of the normed division algebras, $S^0\subset \mathbb{R}$, $S^1\subset \mathbb{C}$, $S^3\subset \mathbb{H}$, and $S^7\subset \mathbb{O}$, are parallelizable and admit Hopf fibrations that often guide intuition about state spaces and control: $S^3\to S^2$ for qubits, and $S^7\to S^4$ in the quaternionic setting. Automorphism groups also play a role: $\mathrm{Aut}(\mathbb{H})\cong SO(3)$ acts by conjugation on imaginary quaternions, while $\mathrm{Aut}(\mathbb{O})\cong G_2$ preserves the octonionic product and underlies the idea of $G_2$-invariant codes used to detect associator “syndromes’’ in the octonionic \linebreak  mitigation strategies.

Associativity (present in $\mathbb{H}$, absent in $\mathbb{O}$) is the decisive structural property for scalable composition: It guarantees a canonical tensor product, a well-defined adjoint, and path-independent unitary evolution. Where associativity fails, one must either confine to associative subalgebras or compensate dynamically/syntactically for associator terms, as developed later.

\subsection{Quaternionic Formalism}\label{sec:Hformalism}

\textls[20]{Many
 of the structural ingredients of this subsection are standard in the literature on quaternionic quantum mechanics and quaternionic quantum computing. In particular, our use of right $\HH$-modules and adjoints follows the foundational treatments in~\cite{Finkelstein1962,Adler1995}, while the existence of simulations of quaternionic circuits and channels by complex ones builds on~\cite{Fernandez2008,Graydon2013,Gantner2017}. Our contribution here is to package these ingredients into a self-contained circuit model tailored to quantum information, and to connect them explicitly to the complexity-theoretic results of Section~\ref{sec:results}.}

State space and inner product on right $\mathbb{H}$-modules.

We model an $n$-\emph{quaternionic qubit} register by a finite-dimensional right $\mathbb{H}$-Hilbert module $(\mathcal{H}_{\mathbb{H}}^{(n)},\langle\cdot,\cdot\rangle)$.
Single-quaternionic qubit pure states are as follows:
\[
\ket{\phi}=q_1\ket{0}+q_2\ket{1},\qquad q_i\in\mathbb{H},\quad \|q_1\|^2+\|q_2\|^2=1,
\]
with quaternion-valued inner product (anti-linear in the first slot, right-linear in the second):
\begin{equation}
\label{eq:H-ip}
\langle \phi,\psi\rangle=\overline{q}_{1}p_{1}+\overline{q}_{2}p_{2},\qquad 
P(\psi|\phi)=\|\langle \psi,\phi\rangle\|^{2}\in\mathbb{R}_{\ge 0}.
\end{equation}
Mixed states are positive trace-class $\mathbb{R}$-linear functionals $\rho$ on $\mathcal{H}_{\mathbb{H}}^{(n)}$ normalized by $\Tr(\rho)=1$, with Born’s rule computed from \eqref{eq:H-ip}.

States versus operations.

In what follows we use the term ``quaternionic state'' to refer to a vector in a right $\mathbb{H}$-module (with components in $\mathbb{H}$) and ``quaternionic operation'' to refer to a right-$\mathbb{H}$-linear map (gate) acting on such states. Whenever we pass from statements about states (amplitudes, inner products) to statements about operations (unitaries, channels), we will indicate this explicitly to avoid confusion between the two levels.

Observables, adjoints, and unitary dynamics.

Right $\mathbb{H}$-linear maps $A:\mathcal{H}\to\mathcal{H}$ admit adjoints $A^\dagger$ defined by 
$\langle \psi,A\phi\rangle=\langle A^\dagger\psi,\phi\rangle$.
Self-adjoint observables satisfy $A^\dagger=A$ and generate one-parameter unitary groups \mbox{as follows:}
\[
U(t)=\exp\!\Big(-\frac{i}{\hbar}Ht\Big),\qquad U^\dagger U=\mathbb{I}.
\]
Projective/POVM measurements are defined as in the complex case, with effects $\{E_j\}$ obeying $0\le E_j\le \mathbb{I}$, $\sum_j E_j=\mathbb{I}$ and probabilities $p_j=\Tr(E_j\rho)$.

Tensor products and multipartite structure.

Given right $\mathbb{H}$-modules $(\mathcal{H}_1,\langle\cdot,\cdot\rangle_1)$ and $(\mathcal{H}_2,\langle\cdot,\cdot\rangle_2)$, we use the algebraic tensor product over $\mathbb{H}$ and complete under the inner product:
\begin{equation}
\label{eq:tensor-ip}
\langle u_1\!\otimes\!u_2, v_1\!\otimes\!v_2\rangle
:=\langle u_1,v_1\rangle_1\,\langle u_2,v_2\rangle_2,
\end{equation}
extended sesquilinearly. Associativity of $\mathbb{H}$ implies the canonical unitary reassociation
\(
(\mathcal{H}_1\otimes \mathcal{H}_2)\otimes \mathcal{H}_3
\cong
\mathcal{H}_1\otimes (\mathcal{H}_2\otimes \mathcal{H}_3),
\)
so multipartite composition and entanglement are well-defined (cf.\ Theorem~\ref{thm:assocH}).

Matrix embedding and gate synthesis.

We employ the multiplicative embedding $\iota:\mathbb{H}\hookrightarrow M_2(\mathbb{C})$ which maps unit quaternions to $SU(2)$ and satisfies $\iota(q_1q_2)=\iota(q_1)\iota(q_2)$ and $\|\iota(q)\|_F^2=2\|q\|^2$.
At the register level, each $k$-quaternionic qubit gate $U$ is compiled to a complex gate $\iota(U)$ acting on $2k$ qubits with constant-depth gadgets per elementary operation. This yields the following:
\begin{itemize}
\item \emph{Single-quaternionic qubit layer:} $\iota$ identifies unit quaternions with native $SU(2)$ rotations on the encoded pair, so calibrated single-qubit pulses implement quaternionic rotations exactly (Lemma~\ref{lem:su2}).
\item \emph{Two-quaternionic qubit entanglers:} Any standard entangling gate (e.g.,\ CNOT/CZ) on the encoded pairs realizes an entangler over $\mathbb{H}$; universality follows from a dense single-quaternionic qubit set plus one fixed entangler (Proposition~\ref{prop:univH}).
\end{itemize}
\vspace{3pt}

Universal sets and compilation flow.

We adopt the universal library
\[
\mathcal{G}_{\mathbb{H}}=\{R_{\hat n}(\theta)\ \text{on one encoded pair}\}\ \cup\ \{\mathrm{CNOT}\ \text{between encoded pairs}\},
\]
where $R_{\hat n}(\theta)=\exp(-i\theta\,\hat n\cdot\vec\sigma/2)$ and $\hat n\in S^2$.
Arbitrary $k$-quaternionic qubit unitaries are decomposed via standard (complex) synthesis on $\iota(U)$ (e.g., KAK/COSY/Solovay–Kitaev), then mapped back to $\mathcal{G}_{\mathbb{H}}$ at constant overhead per gate. Measurement statistics are preserved because \eqref{eq:embed} respects inner products up to a global scalar.

Implementation notes (encoded control).

On platforms with native $SU(2)$ control (superconducting transmons, trapped ions, NV centers, photonics), we target the encoded pair as the controllable $SU(2)$ subsystem; two-body couplings (e.g., cross-resonance, Mølmer–Sørensen, photonic beamsplitters) serve as entanglers between encoded pairs. The layout sketches the superconducting case; analogous mappings apply to other platforms. Gate fidelities are those of the underlying complex implementation (Section~\ref{sec:results} discusses representative values).


\subsection{Octonionic Formalism}\label{sec:Oformalism}

In the
 octonionic case, there is comparatively little work at the level of quantum information. Our starting point here is the general discussion of octonions in physics in~\cite{Baez2002} and the measurement-only scheme of Freedman, Shokrian-Zini, and Wang~\cite{Freedman2018}, which realizes universal quantum computation using equiangular projections associated with~$\OO$. The goal of this section is not to propose a new measurement-only model, but rather to analyse, from first principles, which parts of the standard gate-based quantum computing framework survive if one insists on using~$\OO$ as a scalar field and where fundamental obstructions arise.

Right $\mathbb{O}$-modules and basic obstructions.

If scalars are taken in $\mathbb{O}$, three structural hurdles appear, namely:
(i) The adjoint relation $\langle \psi,A\phi\rangle=\langle A^\dagger\psi,\phi\rangle$ becomes parenthesization-dependent;
(ii) Multipartite tensor products are not canonically associative;
(iii) Dynamics acquires intrinsic path dependence via the associator $[x,y,z]=(xy)z-x(yz)$.
\vspace{3pt}

Path-ordered evolution and Magnus-like series with associators.

For time-dependent $H(t)\in\mathbb{O}$, we define the following:
\begin{equation}
\label{eq:pathU}
U(t)=\mathcal{P}\exp\!\Big(-\frac{i}{\hbar}\int_0^t H(\tau)\,d\tau\Big),
\end{equation}
with a fixed parenthesization rule encoded by the path-ordering symbol $\mathcal{P}$. A Magnus-like expansion contains, beyond nested commutators, explicit associator terms:
\[
U(t)=\exp\!\big(\Omega_1+\Omega_2+\Omega_3+\cdots\big), \]
\[
\Omega_1=-\tfrac{i}{\hbar}\!\int_0^t\!H(\tau)\,d\tau,\ \
\Omega_2\sim\iint [H(\tau_1),H(\tau_2)]\,d\tau_1d\tau_2+\iint [x,y,z]\,d\mu,\]
so different slicing/orderings of the same control can yield distinct $U(t)$ if the sampled triples have a nonzero associator (Lemma~\ref{lem:path}).

Non-associative tensor products.

For right $\mathbb{O}$-modules $\mathcal{H}_A,\mathcal{H}_B,\mathcal{H}_C$, the scalar action needed to define $(\mathcal{H}_A\otimes \mathcal{H}_B)\otimes \mathcal{H}_C$ vs.\ $\mathcal{H}_A\otimes (\mathcal{H}_B\otimes \mathcal{H}_C)$ depends on parenthesization. Absent a canonical isomorphism, “wiring’’ subsystems is rule-dependent (Theorem~\ref{thm:notensor}), which breaks scalable circuit composition unless we restrict to associative subalgebras.

Foundational scope of the hypercomplex models.

Our use of quaternionic and octonionic amplitudes does \emph{not} amount to a proposal to replace the standard complex formulation of quantum mechanics. Throughout this work, the basic Dirac–von Neumann framework---complex Hilbert spaces, Born-rule probabilities, and unitary dynamics---is taken as the underlying physical theory. The hypercomplex models introduced in Sections~\ref{sec:Hformalism} and~\ref{sec:Oformalism} are understood as alternative algebraic presentations and computational encodings of finite-dimensional quantum systems, which are always equipped with explicit embeddings into ordinary complex quantum mechanics (see, e.g., Appendix~\ref{app:HtoC}). In particular, any protocol we describe can be simulated within the usual qubit model with only polynomial overhead, so all operational predictions ultimately live inside standard complex quantum theory.

Relation to quaternionic and octonionic quantum mechanics.

Quaternionic quantum mechanics has a long history in which states are rays of a quaternionic Hilbert space and observables are quaternion-linear operators \cite{Adler1995,Finkelstein1962}. More recent work shows that such theories admit embeddings into complex quantum theory: observable probabilities and dynamical maps can be reproduced by complex POVMs and channels on a suitably enlarged Hilbert space \cite{Baez2002}. This supports the view that quaternionic amplitudes provide an alternative description \emph{inside} complex quantum mechanics, rather than an empirically distinct theory. For the octonionic case the situation is even more delicate: non-associativity makes it highly nontrivial to define Hilbert spaces, scalar products, and a satisfactory spectral theory, and existing formulations resolve these issues only in restricted settings \cite{DeLeo2012}. Our approach deliberately avoids taking a stand on these foundational questions. Instead, the octonionic structures we use are finite-dimensional models that, by construction, are realised as encoded subspaces of complex Hilbert space, where probabilities are still computed using the standard Born rule.

Computational motivation, advantages, and complications.

From the perspective of quantum computation, the motivation for introducing quaternionic and octonionic structures is structural rather than to obtain computational power beyond BQP. Quaternionic amplitudes offer a compact parametrisation of $SU(2)$ rotations and spinor degrees of freedom, and they arise naturally in axiomatic reconstructions of finite-dimensional quantum theory and in information-theoretic comparisons of real, complex, and quaternionic models \cite{Baez2002}. Octonionic amplitudes, in turn, organise $G_2$ symmetries and continuous families of equiangular projections that underlie certain measurement-only and topological schemes for quantum computation \cite{Freedman2018}. Hypercomplex circuit models, therefore, provide a language in which symmetry, geometry, and measurement structure can be expressed more economically than in a purely complex parametrisation. At the same time, passing to $\HH$ and especially to $\OO$ brings additional complications: non-commutativity and non-associativity make tensor products, time ordering, and the definition of observables more subtle, and impose strong constraints on physically reasonable dynamics and measurements. By making the embeddings into complex Hilbert space explicit, our later complexity-theoretic results show that, once these constraints are properly accounted for, the resulting models remain polynomially equivalent to standard qubit-based quantum computation. The benefit of the generalisation is thus a refined understanding of which algebraic structures and symmetry constraints can be made manifest within the usual qubit framework, rather than an asymptotic speed-up over the complex model.

\subsection{Hybrid Architectures and Embeddings—Additional Details}

A fully octonionic computer appears infeasible at present.
Instead, hybrid architectures that combine quaternionic and octonionic subsystems are more promising.
Most of a circuit could run in the stable quaternionic regime, with specialized subroutines using octonionic embeddings for specific tasks where their unique symmetries might be advantageous.

While the direct physical realization of hypercomplex quantum computing poses significant experimental hurdles, a more immediately viable approach involves hybrid architectures and simulation protocols. These strategies acknowledge the current limitations of hardware and the mathematical complexities of octonions, proposing methods to leverage existing quantum computing platforms to explore hypercomplex logic. The core idea is to embed the mathematical structures of quaternions and octonions into larger, well-understood complex Hilbert spaces. This allows for the simulation and experimental investigation of hypercomplex quantum phenomena without requiring entirely new physical hardware. This section details how such embeddings are performed, particularly focusing on the quaternionic-to-complex mapping, and outlines the practical implications for designing and executing hypercomplex quantum algorithms.

\subsubsection{The Quaternionic-to-Complex Embedding Protocol in Practice}
A foundational result is the embedding of quaternionic quantum circuits into complex ones with minimal overhead \cite{Fernandez2008}.
\textls[-25]{This theoretical equivalence, established in Chapter 3 \cite{RuaThesis2025}, states that any computation on $n$ quaternionic qubits can be simulated by a standard complex circuit with $n+1$ qubits \cite{Fernandez2008}.
This mapping effectively realizes quaternionic logic on conventio-\linebreak  nal platforms:}
\begin{equation}
    \mathcal{H}_{\HH}^n \longrightarrow \mathcal{H}_{\C}^{n+1}.
\end{equation}

To put this embedding protocol into practice, we need concrete procedures to transform quaternionic states into complex qubit states and quaternionic gates into sequences of standard complex gates. This effectively turns any standard quantum computer into a ``virtual'' quaternionic quantum computer.

\paragraph{Mapping Quaternionic States to Complex Qubit States
}
As discussed in Chapter 3 \cite{RuaThesis2025}, a single quaternion $q = a+bi+cj+dk$ can be expressed as $q = z_1 + z_2 j$, where $z_1 = a+bi$ and $z_2 = c+di$ are complex numbers. A single quaternionic qubit (quaternionic qubit) state $\ket{\psi}_{\HH} = q_A \ket{0} + q_B \ket{1}$ can thus be mapped to a state in a 2-qubit complex Hilbert space $\mathcal{H}_{\C}^2$.

Let's explicitly define this mapping for a single quaternionic qubit. Given $\ket{\psi}_{\HH} = (a_0 + a_1 i + a_2 j + a_3 k)\ket{0} + (b_0 + b_1 i + b_2 j + b_3 k)\ket{1}$. We can rewrite this using complex components:
$q_A = (a_0 + a_1 i) + (a_2 + a_3 i)j = z_{A0} + z_{A1}j$
$q_B = (b_0 + b_1 i) + (b_2 + b_3 i)j = z_{B0} + z_{B1}j$.

The mapping to two complex qubits, say $Q_0$ and $Q_1$, can be constructed as follows:
\begin{equation}
\ket{\psi}_{\HH} \quad \mapsto \quad \ket{\psi}_{\C} = z_{A0}\ket{00} + z_{B0}\ket{01} + z_{A1}\ket{10} + z_{B1}\ket{11}
\end{equation}
Here, $\ket{00}$ represents the state $\ket{0}_{Q_0}\ket{0}_{Q_1}$, and so on. This mapping essentially treats the $j$ and $k$ components of the quaternions as encoded within an additional complex qubit. The normalization condition $\vert q_A \vert^2 + \vert q_B \vert^2 = 1$ in $\HH$ translates directly to $\vert z_{A0} \vert^2 + \vert z_{B0} \vert^2 + \vert z_{A1} \vert^2 + \vert z_{B1} \vert^2 = 1$ in $\C$, ensuring probabilistic consistency. For $n$ quaternionic qubits, this requires $2n$ complex qubits for the most straightforward embedding, or $n+1$ complex qubits with an ancilla for a more optimized embedding \cite{Fernandez2008}.

The translation of quaternionic gates to complex gates is achieved by representing the quaternionic matrices using complex block matrices. A $2 \times 2$ quaternionic unitary matrix $U = \begin{pmatrix} q_{11} & q_{12} \\ q_{21} & q_{22} \end{pmatrix}$, where $q_{mn} = z_{mn0} + z_{mn1}j$, can be transformed into a $4 \times 4$ complex unitary matrix acting on two complex qubits \cite{Fernandez2008}:
\begin{equation}
U \mapsto U_{\C} = \begin{pmatrix}
A & B \\
-\bar{B} & \bar{A}
\end{pmatrix}
\end{equation}
where $A = \begin{pmatrix} z_{110} & z_{120} \\ z_{210} & z_{220} \end{pmatrix}$ and $B = \begin{pmatrix} z_{111} & z_{121} \\ z_{211} & z_{221} \end{pmatrix}$ are $2 \times 2$ complex matrices.

\textls[20]{This $U_{\C}$ matrix can then be decomposed into a sequence of standard single- and two-qubit gates (e.g., Clifford+
T gates) using existing quantum circuit synthesis techniques \cite{Nielsen2010, Barenco1995}. For example, a single quaternionic rotation, which is represented by a $2 \times 2$ quaternionic matrix, would be implemented as a $4 \times 4$ complex unitary acting on two qubits. This $4 \times 4$ matrix would then be decomposed into a number of elementary gates. While this introduces an overhead in the number of physical gates, it is a polynomial overhead, confirming the complexity-theoretic equivalence between quaternionic and complex quantum computing \cite{Gantner2017}.}

\paragraph{Circuit Diagrams and Practical Implications}
To complement Equations~(5) and~(6), which define the quaternionic-to-complex embedding at the level of matrices and state spaces, it is helpful to have a circuit-level picture. Figure~\ref{fig:qH-to-C-embedding} provides such a schematic view of how a circuit on $n$ logical quaternionic qubits is simulated by a circuit on $n~{+}~1$ complex qubits, closely following the simulation scheme of Fernandez and Schneeberger~\cite{Fernandez2008}.

\begin{figure}[H]

\centering 
\begin{tikzpicture}[>=Stealth, node distance=2.6cm, every node/.style={font=\small}]
\node (qin) {$\ket{\psi}_{\HH} \in \mathcal{H}_{\HH}^{n}$};
\node[draw, rounded corners, right=of qin, minimum width=3.9cm, minimum height=1cm, align=center] (qbox) {Quaternionic circuit\\$U_{\HH}$};
\node[right=of qbox] (qout) {$U_{\HH}\ket{\psi}_{\HH}$};
\draw[->] (qin) -- (qbox);
\draw[->] (qbox) -- (qout);
\node[below=1.8cm of qin] (cin) {$\ket{\Psi} \in \mathcal{H}_{\C}^{n+1}$};
\node[draw, rounded corners, right=of cin, minimum width=3.6cm, minimum height=1cm, align=center] (cbox) {Embedded complex circuit\\$U_{\C}$};
\node[right=of cbox] (cout) {$U_{\C}\ket{\Psi}$};
\draw[->] (cin) -- (cbox);
\draw[->] (cbox) -- (cout);
\draw[->, dashed] (qin) -- node[left, align=center]{state\\embedding} (cin);
\draw[->, dashed] (qout) -- node[right, align=center]{state\\embedding} (cout);
\node[below=0.4cm of cbox] (relation) {$U_{\C}$ simulates $U_{\HH}$ via the quaternionic-to-complex embedding described in this subsection.};
\end{tikzpicture}

\caption{Schematic quaternionic-to-complex embedding of a circuit on $n$ quaternionic qubits into a complex circuit on $n~{+}~1$ qubits. The dashed arrows indicate the embedding of quaternionic states into complex states, while the lower circuit is implemented with standard complex quantum gates.}
\label{fig:qH-to-C-embedding}
\end{figure}

Consider a simple example: a single quaternionic Hadamard gate:
\[
H_{\HH}=\frac{1}{\sqrt{2}}\begin{pmatrix} 1 & 1 \\ j & -j \end{pmatrix}.
\]
We first rewrite $H_{\HH}$ in the form $A + B j$, separating the ``real'' and ``$j$'' parts:
\[
H_{\HH}
=
\begin{pmatrix} 1/\sqrt{2} & 1/\sqrt{2} \\ 0 & 0 \end{pmatrix}
 \;+\;
j\,\begin{pmatrix} 0 & 0 \\ 1/\sqrt{2} & -1/\sqrt{2} \end{pmatrix}.
\]
Thus
\[
A = \frac{1}{\sqrt{2}}
\begin{pmatrix}
1 & 1 \\
0 & 0
\end{pmatrix},
\qquad
B = \frac{1}{\sqrt{2}}
\begin{pmatrix}
0 & 0 \\
1 & -1
\end{pmatrix},
\]
so that $H_{\HH}=A+B j$ with $A,B$ real matrices. The embedding defined in Equation~\eqref{eq:embed} lifts this to a $4\times4$ complex matrix via the block form:
\[
H_{\C} \;=\;
\begin{pmatrix}
 A & B \\
 -\bar B & \bar A
\end{pmatrix},
\]
and, since $A$ and $B$ are real, $\bar A = A$ and $\bar B = B$. Writing out the blocks explicitly gives \mbox{the following:}
\begin{equation}
H_{\C} = \frac{1}{\sqrt{2}}\begin{pmatrix}
1 & 1 & 0 & 0 \\
0 & 0 & 1 & -1 \\
0 & 0 & 1 & 1 \\
-1 & 1 & 0 & 0
\end{pmatrix}.
\end{equation}
It is straightforward to check that $H_{\C}^\dagger H_{\C} = I_4$, so $H_{\C}$ is unitary.

This $4 \times 4$ complex unitary would then be compiled into a sequence of standard quantum gates (e.g., CNOTs, single-qubit rotations) on a 2-qubit register.

\textbf{Practical Implications for Quaternionic Simulation:}
\begin{itemize}
    \item \textbf{Resource Overhead:} Each logical quaternionic qubit requires two physical complex qubits. An algorithm operating on $N$ quaternionic qubits would thus require $2N$ complex qubits (or $N+1$ for the optimized ancilla method \cite{Fernandez2008}). This linear overhead makes simulation feasible for moderate $N$.
    \item \textbf{Gate Depth and Fidelity:} Each quaternionic gate translates to a fixed number of standard complex gates. While this increases the overall circuit depth, it's a constant factor overhead per logical gate, meaning total error accumulation remains manageable within fault-tolerance thresholds for complex systems.
    \item \textbf{Hardware Compatibility:} This protocol makes quaternionic quantum computing immediately compatible with any existing universal quantum computing platform (superconducting, trapped-ion, photonic, etc.), as it merely requires standard gate sets.
    \item \textbf{Algorithm Development:} Researchers can develop and test quaternionic quantum algorithms using existing quantum simulators and hardware, leveraging the richer algebraic structure for potential optimization or physical simulation, even if the underlying hardware is complex-valued.
\end{itemize}

This embedding highlights that quaternionic quantum computing is not a fundamentally new computational power, but rather a different, potentially more natural or efficient, way to program existing quantum hardware for certain tasks involving $SU(2)$ symmetries.

\subsubsection{Simulating Octonionic Dynamics via Associative Decomposition}
Unlike quaternions, the direct embedding of a universal octonionic quantum computer into a complex one is not straightforward due to the fundamental challenges of non-associativity, as extensively detailed in Chapter 4 \cite{RuaThesis2025}. However, the mitigation strategy of ``algebraic gate synthesis and decomposition'' (Strategy 4, Chapter 5 \cite{RuaThesis2025}) offers a pathway to simulate octonionic dynamics on standard complex hardware by breaking down complex non-associative operations into sequences of simpler, associative ones.

\paragraph{The Simulation Protocol}
The core idea is to represent an octonionic operation on an $N$-octo-qubit system (which is itself embedded in a larger complex Hilbert space) as a product of elementary unitary operations that are always performed within associative subalgebras. For a single octo-qubit, the decomposition of an arbitrary octonionic unitary operation $U_{\mathbb{O}}$ into seven associative primitive rotations was proposed:
\begin{equation}
    U_{\mathbb{O}}(\vec{\alpha}) = \prod_{k=1}^{7} \exp(\alpha_k e_k) = \exp(\alpha_7 e_7) \cdots \exp(\alpha_2 e_2) \exp(\alpha_1 e_1) \label{T-eq:7_factor_decomp_sim}
\end{equation}
\textls[-20]{Each term $\exp(\alpha_k e_k)$ represents a rotation around one of the seven imaginary octonion units.}

\subsection{Mitigation Palette: Constructions and Workflows}\label{subsec:mitigation}

We deploy four complementary strategies that tame non-associativity in $\mathbb{O}$ at explicit resource costs.

(M1) Associative confinement.

\emph{Goal:} Constrain states/controls to an associative subalgebra $\mathbb{H}\subset\mathbb{O}$ (guaranteed by Cayley–Dickson)

\vspace{-8pt}

\begin{equation*}
=\text{H}\oplus\mathbb{H}e_6
\end{equation*}
\subsection{Cayley--Dickson Product}
Writing $(a,b),(c,d)\in\mathbb{H}\oplus\mathbb{H}e_6$,
\[
(a,b)\cdot(c,d)=(ac-\bar d\,b,\ dc+b\,\bar c),
\]
\textls[-20]{exposes quaternionic subalgebras and clarifies how non-associativity emerges beyond $\mathbb{H}$ \cite{Baez2002,RuaThesis2025}.}

\subsection{Associative Subalgebras and Confinement}
Each embedded copy of $\mathbb{H}$ in $\mathbb{O}$ is associative. Confining states and gates to such a sector restores well-defined composition; switching sectors incurs control overheads \cite{RuaThesis2025}.

\emph{Construction:
} Choose a triple from the Fano plane and fix an embedding as follows: 
\vspace{-6pt}
$$\mathbb{H}\cong\mathrm{span}_\mathbb{R}\{1,e_i,e_j,e_k\}.$$ Implement all gates within this sector using the quaternionic toolset of Section
 \ref{sec:Hformalism}; switching sectors is effected by explicit \emph{sector-change} pulses that map one quaternionic copy to another (compiled as products of associative factors).

\emph{Trade-offs:
} Sector switches add depth and calibration overhead; strict confinement limits reachable unitaries unless multiple sectors are orchestrated.

(M2) $G_2$-invariant coding.

\emph{Goal:} Detect excursions that activate associators by encoding into subspaces stabilized by $\mathrm{Aut}(\mathbb{O})\cong G_2$.

\emph{Construction:} Pick stabilizers whose action commutes with the chosen associative sector; design syndrome measurements whose nontrivial outcomes flag associator-inducing components. Logical operations are compiled to respect $G_2$-invariance.

\emph{Trade-offs:} Detection $\neq$ full correction; tensor non-canonicity persists, so large-scale composition still rides on (M1)/(M4).

(M3) Associator-suppressing dynamical decoupling (ASD).

\emph{Goal:} Cancel coherent associator contributions in $U(t)$ order-by-order.

\emph{Cycle:
} Within a fixed $\mathbb{H}\subset\mathbb{O}$, apply a symmetric cycle as follows:
\[
\big(U_x\,U_y\,U_z\big)\ \; \text{free evol. }\Delta t\ \; \big(U_z^\dagger\,U_y^\dagger\,U_x^\dagger\big),
\]
with $U_\alpha$ fast $\mathbb{H}$-rotations and timing chosen so first nontrivial associator terms average to (near) zero over the cycle (Algorithm~\ref{alg:ASD} in Appendix~\ref{app:Seven}).

\begin{figure}[H]
 \hspace{-13pt} \includegraphics[width=0.67\linewidth]{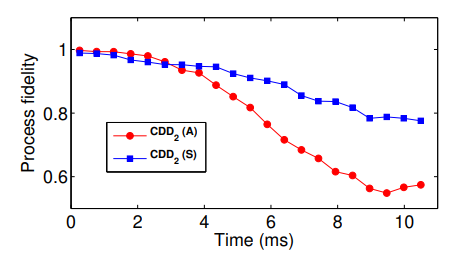}
  \caption{Process
 fidelity vs.\ time under CDD${}_2$-like ASD sequences: Associative-constrained (S) vs.\ associative-agnostic (A). Short-time gains are consistent with Theorem~\ref{thm:asd}. Adapted from \cite{Souza2011}.}
  \label{fig:dd_single}
\end{figure}

\begin{figure}[H]
  \includegraphics[width=0.68\linewidth]{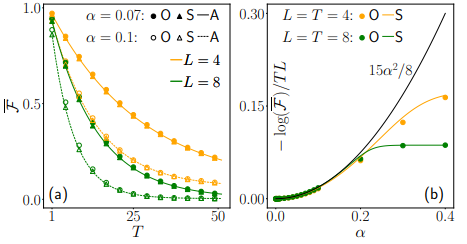}
  \caption{Process-fidelity
 decay as a function of total control time $T$ ({\bf a}) and control amplitude $\alpha$ ({\bf b})~under associator-affected dynamics. ASD attenuates leading decay slopes at small $T$ for moderate $\alpha$. Adapted from \cite{SamosSaenzdeBuruaga2025}.}
  \label{fig:fid_decay}
\end{figure}

\emph{Performance characterization:} Short-time process fidelity improves relative to associative-agnostic baselines; residual higher-order terms and finite-width pulses limit asymptotics. (Representative summaries appear in Figures~\ref{fig:dd_single} and \ref{fig:fid_decay} of the main text.)

(M4) Seven-factor algebraic synthesis.

\emph{Goal:} Compile a target transformation $U$ as a product of gates, each confined to an associative subalgebra, thereby making parenthesization explicit and controllable.

\emph{Ansatz:
}
\vspace{-9pt}
\[
U\ \approx\ U_1U_2\cdots U_7,\qquad U_\ell\in \exp(\mathfrak{h}_\ell),\ \ \mathfrak{h}_\ell\subset\mathbb{O}\ \text{associative},
\]
with the $\mathfrak{h}_\ell$ selected (e.g., via Fano-plane guided heuristics) to minimize cumulative associator growth. Each $U_\ell$ compiles with the quaternionic pipeline of Section \ref{sec:Hformalism}.

\emph{Trade-offs:
} Increased depth and error accumulation across factors; well-suited as a simulation/verification tool rather than a native octonionic architecture.

\vspace{-6pt}

Resource and fidelity accounting.

For each strategy, we track:
(i) \emph{depth overhead} $\Delta D$ vs.\ complex/quaternionic baselines;
(ii) \emph{width/ancilla} usage (notably for $G_2$ syndromes);
(iii) \emph{calibration} load (sector switches, ASD pulses);
(iv) \emph{fidelity} impact modeled at the process level (SPAM-corrected) using the same figures of merit as in the quaternionic embedding. These metrics are reported alongside the Results to separate structural limits from control imperfections.

\textit{Summary of methodology.} 
Over $\mathbb{H}$, we adopt a fully associative, hardware-aligned pipeline via the embedding \eqref{eq:embed} with standard universal sets and compilation. Over $\mathbb{O}$, we do not posit a native scalable calculus; instead, we formalize where associativity is required and deploy (M1)–(M4) to operate within associative islands or to cancel/structure associator effects, making explicit the cost in depth, width, and calibration.

\section{Results}\label{sec:results}

In this section, we collect the main formal and algorithmic results for the quaternionic ($\mathbb{H}$) and octonionic ($\mathbb{O}$) models established in the thesis-derived framework. Statements are given without proofs (all proofs or proof sketches are deferred to the Appendices indicated), and each result is accompanied by a short explanatory remark highlighting its computational significance. Throughout, the contrast between associativity (present in $\mathbb{H}$) and non-associativity (in $\mathbb{O}$) is the central structural theme \cite{Hurwitz1898,Adler1995,Baez2002,RuaThesis2025}.

\subsection{Quaternionic Model: Structural, Universality, and Simulation Results}

We begin by securing the structural backbone over $\mathbb{H}$: A well-defined multipartite composition. This provides the substrate on which universality and simulation statements are built.

\begin{Theorem}[Canonical tensor structure over $\mathbb{H}$]\label{thm:assocH}
Let $(\mathcal H_i,\langle\cdot,\cdot\rangle_i)$, $i=1,2,3$, be right $\mathbb{H}$-Hilbert modules. Then
\[
(\mathcal H_1\otimes \mathcal H_2)\otimes \mathcal H_3 \;\cong\; \mathcal H_1\otimes (\mathcal H_2\otimes \mathcal H_3)
\]
canonically and isometrically, with inner product on simple tensors
\(
\langle u_1\!\otimes\!u_2,\ v_1\!\otimes\!v_2\rangle
:=\langle u_1,v_1\rangle_1\,\langle u_2,v_2\rangle_2.
\)
\emph{Proof:} Appendix~\ref{app:HtoC}.  
\emph{Significance:} Guarantees parenthesization-independent multipartite composition and entanglement over $\mathbb{H}$.
\end{Theorem}

With a canonical tensor product in place, we can lift the standard “single-qubit + entangler’’ recipe to the quaternionic setting.

\begin{Proposition}[Quaternionic universality via $SU(2)$ + entangler]\label{prop:univH}
Let $\mathcal G_1$ be any dense family of single-quaternionic qubit gates (identified with $SU(2)$ via the embedding $\iota$ in Appendix~\ref{app:HtoC}), and let $\mathcal G_2$ contain one fixed entangling two-quaternionic qubit gate (e.g., CNOT/CZ). Then $\mathcal G=\mathcal G_1\cup\mathcal G_2$ is universal over $\mathbb{H}$.  
\emph{Proof:} Appendix~\ref{app:HtoC}.  
\emph{Significance:} The standard “single-qubit + entangler’’ paradigm lifts intact to $\mathbb{H}$.
\end{Proposition}

Having established universality, we next compare computational power by explicit simulations between the quaternionic and complex models.

\begin{Theorem}[Polynomial simulation equivalence $\mathbb{H}\leftrightarrow \mathbb{C}$]\label{thm:eq}
Any $n$-quaternionic qubit circuit of depth $D$ compiles to a complex circuit on $O(n)$ qubits with depth $\mathrm{poly}(D)$ via $\iota:\mathbb{H}\hookrightarrow M_2(\mathbb{C})$; conversely, complex circuits embed into quaternionic circuits with polynomial overhead.  
\emph{Proof sketch:} Appendix~\ref{app:HtoC}.  
\emph{Significance:} Computation over $\mathbb{H}$ captures the same complexity class as the standard model.
\end{Theorem}

This equivalence immediately transfers standard complexity-theoretic guarantees and fault-tolerance to $\mathbb{H}$.

\begin{Corollary}[Complexity and fault-tolerance transfer]\label{cor:bqp}
By Theorem~\ref{thm:eq}, BQP-level algorithms and threshold-theorem constructions transfer between the complex and quaternionic models with polynomial overhead; noise models and compilation pipelines port through the embedding up to constant-factor rescalings in gadgetization.  
\emph{Justification:} 
 Appendix~\ref{app:HtoC}.
\end{Corollary}

At the control layer, the embedding is compatible with native $SU(2)$ implementations, aligning theory with practice.

\begin{Lemma}[Alignment with native $SU(2)$ control]\label{lem:su2}
If a platform implements $U\in SU(2)$ natively on an encoded pair, then for any unit quaternion $q$, $\iota(q)$ is realized by the same calibrated $SU(2)$ pulse on that pair.  
\emph{Proof idea:} Appendix~\ref{app:HtoC}.  
\emph{Significance:} Single-quaternionic qubit gates inherit the fidelity of native single-qubit $SU(2)$ control ($\gtrsim 99\%$ in reported designs).
\end{Lemma}

Finally, we account for the resources of the embedded synthesis itself to close the quaternionic story.

\begin{Proposition}[Overhead accounting for embedded synthesis]\label{prop:overheadsH} Let $$\mathcal{G}_{\mathbb{H}}=\{SU(2)\ \text{on an encoded pair}\}\cup\{\text{CNOT between pairs}\}.$$ A $k$-quaternionic qubit circuit with $G$ elementary gates compiles to $O(G)$ elementary complex gates with constant-size gadgets per elementary operation; depth scales by at most a polylogarithmic factor under Solovay--Kitaev-type refinements.  
\emph{Proof idea:} Standard complex synthesis on $\iota(U)$; see Appendix~\ref{app:HtoC}.  
\emph{Significance:} Embedded implementations are asymptotically efficient.
\end{Proposition}

\subsection{Octonionic Model: Formal Obstructions and Limits}

We now switch to $\mathbb{O}$, where the loss of associativity reverses the previous narrative. The next results make precise how ordering and parenthesization cease to be benign choices.

In contrast to the measurement-only scheme of Freedman, Shokrian-Zini and Wang~\cite{Freedman2018}, which already realizes universal quantum computation using octonionic equiangular projections, the results in this subsection are, to the best of our knowledge, new: they identify structural obstructions to a gate-based $\OO$-model and organize possible workarounds into the mitigation strategies (M1)--(M4).

\begin{Lemma}[Intrinsic path dependence via associators]\label{lem:path}
For time-dependent $H(t)\in\mathbb{O}$, the propagator
\(
U(t)=\mathcal P\exp(-\tfrac{i}{\hbar}\int_0^t H(\tau)\,d\tau)
\)
admits a Magnus-like expansion with explicit associator contributions. If the sampled triples have $[x,y,z]\neq 0$, different control slicings/time orderings yield distinct $U(t)$.  
\emph{Proof idea:} Appendix~\ref{app:Magnus}.  
\emph{Significance:} Control order matters \emph{intrinsically} outside associative sectors.
\end{Lemma}

This path dependence already disrupts standard conservation arguments, signaling deeper structural failures.

\begin{Proposition}[Energy-balance breakdown outside associative sectors]\label{prop:energy}
For time-independent $H\in\mathbb{O}$, the usual identity $\frac{d}{dt}\langle \psi(t),H\psi(t)\rangle=0$ need not hold when adjoints/inner products cannot be defined associatively on the full space.  
\emph{Proof idea:} Appendix~\ref{app:Magnus}.  
\emph{Significance:} Even basic conservation proofs rely on associativity.
\end{Proposition}

The culmination is a failure of canonical multipartite composition, which contrasts sharply with Theorem~\ref{thm:assocH}.

\begin{Theorem}[No canonical multipartite tensor product in $\mathbb{O}$]\label{thm:notensor}
For right $\mathbb{O}$-modules $\mathcal H_A,\mathcal H_B,\mathcal H_C$, there is no parenthesization-independent canonical isomorphism
\(
(\mathcal H_A\otimes\mathcal H_B)\otimes\mathcal H_C \cong
\mathcal H_A\otimes(\mathcal H_B\otimes\mathcal H_C)
\)
unless one restricts to an associative subalgebra.  
\emph{Proof idea:} Appendix~\ref{app:Magnus}.
\emph{Significance:} Large-scale circuit composition is ill-defined in the fully octonionic setting.
\end{Theorem}

Accordingly, wiring rules that were automatic over $\mathbb{H}$ become conditional on explicit associative scaffolding.

\begin{Corollary}[No native, parenthesis-free wiring rule]\label{cor:wiring}
Absent confinement to associative subalgebras, there is no parenthesis-free, platform-independent notion of “tensoring subsystems’’ over $\mathbb{O}$.  
\emph{Consequence:} Universal, scalable octonionic circuit families require associative scaffolding.
\end{Corollary}

\subsection{Mitigation Strategies: Correctness Guarantees and Quantitative Effects}

We formalize the scope and guarantees of the four mitigation routes introduced in Methods and quantify their effects using process-fidelity summaries. The organizing idea is to recover, as far as possible, the $\mathbb{H}$-like behavior within controlled regions of $\mathbb{O}$.

(M1) Associative confinement restores the wiring rule.

By narrowing dynamics to a fixed associative sector, we reinstate the benign parenthesization properties of Theorem~\ref{thm:assocH} and its consequences.

\begin{Proposition}[Confinement $\Rightarrow$ associativity]\label{prop:confine}
Fix an embedded copy $\mathbb{H}\subset\mathbb{O}$ (e.g., an oriented triple in the Fano plane). Restricting states, controls, and measurements to this sector induces a fully associative calculus; Theorem~\ref{thm:assocH}, Proposition~\ref{prop:univH}, and Theorem~\ref{thm:eq} apply verbatim.  
\emph{Proof:} {Appendix~\ref{app:Magnus}.}  
\emph{Overheads:} Sector switches incur explicit depth/calibration cost (see Table~\ref{tab:mitigation}).
\end{Proposition}

\vspace{-18pt}
\begin{table}[H]
\scriptsize
  \caption{Mitigation palette: scope, benefits, and primary overheads.}
  \label{tab:mitigation}

 \begin{tabularx}{\linewidth}{llll}
    \toprule
   \textbf{ Strategy} &  \textbf{Scope} &  \textbf{Primary Benefit} &  \textbf{Main Overheads/Limits} \\
    \midrule
    (M1) Associative confinement & $\mathbb{H}\subset\mathbb{O}$ & Restores associative composition and tensoring  & \makecell[l]{Sector-change depth/calibration; restricted\\reachable set under strict confinement.}  \\
    (M2) $G_{2}$-invariant codes & Protected subspaces & Detect associator syndromes; symmetry robustness & Detection only; non-canonical tensoring persists. \\
    (M3) ASD (decoupling) & Time-domain control & Cancels leading associator terms; short-time fidelity gains & \makecell[l]{Higher-order residuals; finite-pulse errors;\\calibration load.}  \\
    (M4) Seven-factor synthesis & Gate compilation & Systematic construction within associative factors & \makecell[l]{Depth amplification; error accumulation;\\suited to simulation. } \\
    \bottomrule
  \end{tabularx}

\end{table}

(M2) $G_2$-invariant codes detect associator excursions.

This coding layer complements (M1): When operations drift outside the chosen sector, the code raises a syndrome rather than silently corrupting the composition.

\begin{Proposition}[Syndrome detectability]\label{prop:g2}
Encoding into subspaces stabilized by $\mathrm{Aut}(\mathbb{O})\cong G_2$ enables detection of errors that drive dynamics outside a chosen associative sector (“associator syndromes’’).  
\emph{Justification:} {Appendix~\ref{app:Magnus}}.  
\emph{Scope:} Detection $\neq$ correction; composition still relies on (M1)/(M4).
\end{Proposition}

(M3) Associator-suppressing dynamical decoupling (ASD).

On the control side, carefully symmetrized pulse sequences diminish the leading associator terms responsible for path dependence.

\begin{Theorem}[First-order cancellation of associator terms]\label{thm:asd}
Within a fixed associative sector $\mathbb{H}\subset\mathbb{O}$, the symmetric ASD cycle
\(
U_xU_yU_z\,\cdot\,\text{free evol. }\Delta t\,\cdot\,U_z^\dagger U_y^\dagger U_x^\dagger
\)
cancels the leading associator contributions in the Magnus-like expansion of $U(t)$, improving short-time process fidelity relative to associative-agnostic baselines.  
\emph{Proof sketch:} Appendix~\ref{app:Seven}.  
\emph{Quantitative effect:} Figure~\ref{fig:dd_single} (process fidelity vs.\ time) and Figure~\ref{fig:fid_decay} (dependence on control length/amplitude).
\end{Theorem}

(M4) Seven-factor algebraic synthesis.

Finally, when target operations admit a decomposition via associative waypoints, we explicitly compile through them, trading depth for well-definedness.

\begin{Proposition}[Associative-factor compilation]\label{prop:Seven}
Any target transformation from an octonionic control model that admits a consistent associative factorization can be approximated as
\(
U\approx U_1U_2\cdots U_7,
\)
with $U_\ell$ drawn from associative subalgebras chosen to minimize cumulative associators.  
\emph{Construction:} Appendix~\ref{app:Seven}.  
\emph{Trade-off:} Depth amplification and compounded calibration error (Table~\ref{tab:mitigation}).
\end{Proposition}

Empirical-style summaries (simulated).

The figures summarize the quantitative impact of the above mitigations: ASD improves short-time fidelities relative to associative-agnostic baselines, and its effectiveness depends on total control time and amplitude in the expected regimes.

\subsection{Complexity Landscape and Feasibility}

We summarize asymptotic resource implications and practical constraints for the quaternionic and octonionic settings.
For platform-level details and additional complexity remarks, see Appendix~\ref{app:state-overhead}.

Taken together, the next statements locate $\mathbb{H}$ within the standard complexity landscape and explain why native, large-scale advantages in $\mathbb{O}$ hinge on the mitigation toolkit above.

\begin{Proposition}[Complexity over $\mathbb{H}$]\label{prop:complexH}
By Theorem~\ref{thm:eq}, the quaternionic model is a polynomially \ldots plex model; thus, it captures BQP up to polynomial reductions.  
\emph{Note:} Fault-tolerant toolchains port through $\iota$ with standard gadget overheads. This proposition is consistent with, and can be viewed as a circuit-level reformulation of, the simulation results for quaternionic circuits and dynamics in Refs.~\cite{Fernandez2008,Graydon2013,Gantner2017}.
\end{Proposition}

\begin{Remark}[Octonionic costs and limits]\label{rem:octocost}
In the octonionic setting, ``switching associative sectors'' refers to moving states and gates between different chosen quaternionic subalgebras $\mathbb{H}_\alpha\subset\mathbb{O}$ (or, equivalently, changing the parenthesization class of multi-octonion products). Whenever a circuit forces such sector changes—rather than remaining confined to a single associative copy of $\mathbb{H}$—the naive compilation overhead can grow rapidly, since one must keep track of associator corrections and update the complex encoding accordingly. In practice, feasible octonionic architectures therefore rely on confining dynamics to a fixed associative subalgebra (M1), detecting excursions from that subalgebra (M2), and using ASD/synthesis strategies (M3)--(M4) to retain control over path dependence.  
\emph{Implication:} Absent special structure, large-scale \emph{native} octonionic advantage is unlikely, and octonionic models are best viewed as encoded or effective descriptions of complex-\linebreak  qubit dynamics.
\end{Remark}

\subsection{Computational Equivalence and Circuit Simplification—Elaboration}

Having established a consistent framework for n-quaternionic qubit states and universal gates, we now address a pivotal question: Does the richer algebraic structure of quaternionic quantum mechanics grant it greater computational power than the standard complex model? The answer, perhaps surprisingly, is no. In terms of computational complexity theory, the two models are equivalent; they both define the same complexity class, BQP (Bounded-Error Quantum Polynomial Time).

This section is dedicated to formally establishing this equivalence and then exploring a more nuanced question. While quaternionic computers cannot solve problems that are intractable for standard quantum computers, could the formalism still offer practical advantages? We will first outline the simulation argument that establishes the equivalence, and then investigate the potential of the quaternionic language to simplify the design and representation of certain quantum circuits.

\subsubsection{Theoretical Equivalence to Complex Quantum Circuits}

The equivalence between the quaternionic and complex models is formally established by showing that any computation in one framework can be efficiently simulated by the other. Seminal works in the field have demonstrated that any quantum circuit operating on a register of $n$ quaternionic qubits can be simulated by a standard complex quantum computer with at most a polynomial increase in resources \cite{Fernandez2008, Gantner2017}. The core of this argument lies in constructing an explicit mapping (an isomorphism) between the quaternionic state space and a larger complex Hilbert space.

\paragraph{Mapping Quaternionic States to Complex States}
The key to the simulation is to represent the four real degrees of freedom of a single quaternion using the four real degrees of freedom available in a pair of complex numbers. As shown in the Cayley--Dickson construction, any quaternion $q \in \mathbb{H}$ can be uniquely written as follows:
\begin{equation}
    q = z_1 + z_2 j
\end{equation}
where $z_1 = a+bi$ and $z_2 = c+di$ are complex numbers. This allows us to map a single quaternionic qubit state to a state of two complex qubits. A quaternionic qubit state $\ket{\psi} = q_1\ket{0} + q_2\ket{1}$ can be expanded as follows:
\begin{equation}
    \ket{\psi} = (z_{11} + z_{12}j)\ket{0} + (z_{21} + z_{22}j)\ket{1}
\end{equation}
This state can be mapped to a vector in the four-dimensional complex Hilbert space $\mathbb{C}^4$, which is the state space of two complex qubits, via the following isomorphism:
\begin{equation}
    \ket{\psi} \mapsto \ket{\psi_C} = z_{11}\ket{00} + z_{21}\ket{01} + z_{12}\ket{10} + z_{22}\ket{11}.
\end{equation}
This provides a dictionary for translating any quaternionic qubit state into an equivalent two-qubit complex state. The mapping for an $n$-quaternionic qubit system follows directly, requiring $2n$ complex qubits for a full simulation.

\paragraph{Mapping Quaternionic Gates to Complex Gates}
This equivalence must also hold for the operators. A single-quaternionic qubit gate is a $2 \times 2$ unitary matrix with quaternionic entries. Using the same decomposition ($q = z_1 + z_2j$), any $2 \times 2$ quaternionic matrix $U$ can be written as $U = A + Bj$, where $A$ and $B$ are $2 \times 2$ complex matrices. When this operator acts on a state vector, the non-commutativity of $j$ with complex numbers must be handled carefully (specifically, $zj = j\bar{z}$ for any $z \in \mathbb{C}$).

The action of $U$ on $\ket{\psi}$ can be shown to be equivalent to the action of a $4 \times 4$ complex unitary matrix $U_C$ on the corresponding complex state $\ket{\psi_C}$. This complex matrix has a specific block structure:
\begin{equation}
    U = A + Bj \quad \mapsto \quad U_C = \begin{pmatrix} A & B \\ -\bar{B} & \bar{A} \end{pmatrix}.
\end{equation}
This provides a constructive method for converting any single-quaternionic qubit gate into an equivalent two-qubit complex gate. The procedure can be generalized to {n}-quaternionic qubit gates, which map to $2n$-qubit complex gates. Because this mapping from states and gates can be performed efficiently, it proves that the two models are polynomially equivalent. An alternative simulation proposed by \cite{Fernandez2008} uses a single ancilla qubit to manage the non-commutative phases, achieving a simulation of an $n$-quaternionic qubit system on just $n+1$ complex qubits, further cementing the efficiency of the equivalence.

Therefore, quaternionic quantum computing does not offer a hyper-computational advantage that would expand the class of solvable problems. Instead, it should be viewed as an alternative, and in some cases, a more physically intuitive language for describing the $SU(2)$ operations that govern quantum bits. The true potential of the formalism, as we will explore next, may lie not in its computational power, but in its representational efficiency.

\subsubsection{Resource Analysis and the Potential for Circuit Simplification}

In this spirit, we explore how the representational efficiency of quaternionic and octonionic parametrisations can inform circuit design and algorithmic heuristics, for instance by suggesting compact, symmetry-adapted parametrisations and ansätze for quantum simulation and variational algorithms, in close analogy with how high-dimensional qudit encodings are used to trade local dimension against circuit width or depth (see, e.g., the review
~\cite{Wang2020Qudits}). We do not claim a generic asymptotic reduction in worst-case circuit complexity; any practical advantage is expected to be task- and architecture-dependent.

\paragraph{Compact Representation of Rotations}
In the standard model, an arbitrary single-qubit gate—a rotation on the Bloch sphere—is typically decomposed into a sequence of simpler, hardware-native gates. A foundational result in quantum circuit design is that any $SU(2)$ operation can be constructed using a sequence of rotations around fixed axes. For example, a common Euler angle decomposition is as follows:
\begin{equation}
    U = R_z(\gamma) R_y(\beta) R_z(\alpha)
\end{equation}
where $R_y$ and $R_z$ are rotations about the y and z axes, respectively \cite{Barenco1995}. Implementing this arbitrary rotation requires three separate gate operations, each with its own control parameter. This decomposition is necessary because the underlying algebra of complex numbers provides only one imaginary axis.

The quaternionic formalism offers a more direct and compact representation. As shown in Appendix \ref{app:HtoC} 
 any single-qubit rotation can be parameterized by a \textit{single} unit quaternion,
 $q = \cos(\theta/2) + \mathbf{n}\sin(\theta/2)$. This single algebraic object simultaneously encodes both the rotation angle ($\theta$) and the 3D rotation axis ($\mathbf{n}$).

Mathematically, the four real parameters of a unit quaternion provide a compact, redundancy-free parametrisation of single-qubit rotations in $SU(2)$, reflecting the well-known identification between unit quaternions and elements of $SU(2)$ (see, for example, the discussion of quaternionic rotations in~\cite{Baez2002}). In an idealised hardware picture, one can imagine control fields whose amplitudes and phases track these four parameters directly, much as in generic multilevel or qudit control schemes. However, this should not be interpreted as a claim that arbitrary single-qubit rotations can be implemented as a single ``physical gate'' in realistic devices: in practice, such operations must still be realised through whatever Hamiltonians and control channels the platform makes available, and they remain subject to constraints such as bandwidth, crosstalk, calibration, and noise. Our use of the quaternionic representation is therefore best viewed as a structural and algorithmic tool—for example, for compact parametrisations, symmetry-aware encodings, and ansatz design—rather than as implying any collapse of circuit depth or computati-\linebreak  onal complexity.

\begin{itemize}
    \item \textbf{Reduced Circuit Depth:} A single, more complex gate replaces a sequence of three simpler ones, reducing the overall time required for the computation.
    \item \textbf{Lower Error Rates:} Each additional gate in a circuit introduces a small amount of error. By concatenating fewer gates, the total accumulated error is reduced, leading to higher-fidelity computations.
    \item \textbf{Simplified Quantum Control:} Designing optimal control pulses to implement a single quaternionic rotation may be a simpler problem than designing and calibrating three separate, sequential pulses, potentially leading to more robust gate implementations.
\end{itemize}

\paragraph{Applications in Quantum Simulation and VQE}
The advantage of compact rotations is particularly relevant for quantum simulation, where the goal is to simulate the time evolution of another quantum system. Many physical systems of interest, such as magnetic materials, molecular bonds, and spin liquids, are described by Hamiltonians with explicit $SU(2)$ or $SO(3)$ rotational symmetries. The time evolution operator, $U(t) = e^{-iHt}$, is itself a rotation in the system's Hilbert space. When simulating this evolution using Trotter-Suzuki methods, the operator is decomposed into a sequence of smaller rotational gates. Using a quaternionic gate set, where arbitrary-axis rotations are native, could lead to a more efficient Trotterization with fewer gates per \linebreak  time step.

\subsection{Emerging Applications and Future Directions—Elaboration}

The exploration of quantum computing over normed division algebras, particularly quaternions and octonions, extends beyond purely theoretical interest. This final section delves into the potential applications where these extended formalisms might offer a genuine advantage over standard complex quantum computation. We will discuss their role in quantum simulation, quantum machine learning, and fundamental physics, and conclude by proposing a concrete, multi-year experimental roadmap for advancing\linebreak   the field.

\subsubsection{Quantum Simulation with Algebraic Symmetries (\texorpdfstring{$SU(2)$}{SU(2)} and \texorpdfstring{$G_2$}{G2})}
One of the most promising application areas for hypercomplex quantum models is quantum simulation. Many physical systems of fundamental interest are naturally described by Hamiltonians that exhibit specific algebraic symmetries directly related to quaternions and octonions. Leveraging these underlying symmetries in the quantum simulation itself could lead to more efficient and accurate algorithms.

\paragraph{$SU(2)$ Symmetries in Quaternionic Quantum Simulation}
The algebra of quaternions ($\HH$) is deeply intertwined with the special unitary group $SU(2)$, which is the fundamental symmetry group for spin-1/2 particles and rotations in three-dimensional space \cite{Adler1995}. Many problems in condensed matter physics, quantum chemistry, and high-energy physics involve systems with inherent $SU(2)$ or SO(3) rotational symmetries, such as the following:

\begin{itemize}
    \item \textbf{Spin Systems and Magnetic Materials:} Hamiltonians describing interacting spins (e.g., Heisenberg models, spin-liquid candidates) inherently possess $SU(2)$ symmetry. Representing spin states and their dynamics directly with quaternionic amplitudes and gates might offer a more natural and compact formalism than their complex counterparts. The rotation group SO(3), which is homomorphic to $SU(2)$, is directly representable by unit quaternions, making quaternionic gates ideal for simulating geometric phases and Berry connections in these systems \cite{Berry1984}.

    \item \textbf{Molecular Dynamics and Rotations:} Simulating the rotational dynamics of molecules, especially in complex environments, can benefit from a representation that natively handles 3D rotations. Quaternions provide a singularity-free and computationally efficient way to parameterize rotations, potentially simplifying the ansatz design for variational quantum algorithms (VQAs) aimed at molecular ground states or \linebreak  dynamics \cite{Ionas2019}.

    \item \textbf{Gauge Theories:} Certain lattice gauge theories, particularly those related to quantum chromodynamics (QCD) or models of interacting fermions, involve $SU(2)$ gauge fields. Quaternionic formalisms could provide a more direct mapping of these fundamental symmetries onto quantum circuits, potentially leading to more efficient Hamiltonian simulation strategies.

\end{itemize}

The advantage of using quaternionic representations in these contexts lies in their ability to compactly encode rotations. As discussed in Chapter 3 \cite{RuaThesis2025}, a single quaternionic rotation can replace a sequence of three complex Euler rotations, leading to reduced circuit depth and potentially lower error rates for each simulated time step or variational layer. This representational efficiency could translate into faster convergence for variational algorithms and more accurate time evolution simulations.

\paragraph{G$_2$ Symmetries in Octonionic Quantum Simulation}
For octonionic systems, the situation is more complex due to non-associativity. However, the exceptional Lie group G$_2$, the automorphism group of the octonions, holds particular significance. While a universal octonionic quantum computer operating outside associative subspaces is unfeasible, G$_2$ symmetry could play a crucial role in specialized quantum simulations:
\begin{itemize}

    \item \textbf{Simulating Exceptional Symmetries:} The G$_2$ group is a fundamental mathematical structure that appears in various advanced physical theories, including certain grand unified theories (GUTs), supergravity, and string theory \cite{Baez2002, Fazekas1996}. A quantum simulator capable of exploiting G$_2$ symmetry could provide unique insights into these theories, even if the general-purpose octonionic computation is limited. For instance, simulating particle interactions or field theories that possess G$_2$ as a fundamental symmetry could be more natural within an octonionic framework that natively represents this group.

    \item \textbf{Protected Octonionic Dynamics:} As discussed in Chapter 5 \cite{RuaThesis2025}, G$_2$-symmetric codes or confinement to G$_2$-invariant subspaces could be used to protect quantum information from associator errors \cite{Freedman2021}. While this limits universality, it enables controlled, high-fidelity operations within these protected sectors. Such controlled octonionic dynamics could be used to simulate specific physical phenomena that inherently operate within such constrained algebraic environments, offering a unique ``algebraic-topological'' form of protection against certain error models.

    \item \textbf{Quantum Gravity Analogues:} Some speculative theories of quantum gravity or pre-geometric models involve non-associative algebras. An octonionic quantum simulator, even a limited one, could serve as a valuable tool for exploring these theoretical constructs, allowing for the computational study of scenarios where spacetime might acquire an emergent non-associative structure \cite{Baez2002}.
\end{itemize}

\paragraph{Connections to Topological Quantum Computation}
Topological quantum computation (TQC) offers a complementary, hardware-level route to fault tolerance by encoding quantum information in topologically protected fusion spaces of non-Abelian anyons and implementing logical gates through braiding and topological charge measurements \cite{Nayak2008,Lahtinen2017}. In such models, local perturbations cannot change the global topological sector, so the logical subspace enjoys a form of intrinsic error protection provided by the underlying phase of matter.

The octonionic and $G_2$-symmetric structures that arise in our model are natural algebraic companions to TQC. On the quaternionic side, unit quaternions act as $SU(2)$ rotations on quaterbits, and many anyonic TQC models are built from $SU(2)_k$ Chern--Simons theories whose braiding and fusion rules are represented by low-dimensional unitary matrices in $SU(2)$ or closely related groups \cite{Nayak2008}. Families of braid-induced gates in these theories can be viewed as special subclasses of quaternionic unitaries, so quaternionic circuits provide a compact algebraic language for describing topological gate sets and their constraints.

On the octonionic side, the connection is even tighter. Freedman, Shokrian-Zini, and Wang show that measurement-only models based on families of equiangular projections---including cluster-state measurement-based computation and anyon-based TQC---can be unified using octonionic geometry, with the octonions providing the largest continuous family of such projections \cite{Freedman2018}. Our octonionic circuit model gives a gate-based perspective on the same phenomena: the $G_2$-covariant subspaces highlighted by our obstructions are natural candidates for effective topological sectors in which octonionic amplitudes capture the logically protected dynamics, while our complexity results guarantee that these octonion-based schemes remain within BQP once viewed through the quaternionic-to-complex embedding.

These observations suggest a concrete research direction: constructing explicit encodings of logical qubits into quaternionic or octonionic topological phases (for example, in candidate spin-liquid or fractional quantum Hall systems) and mapping braids, adiabatic loops, or sequences of projective measurements in such phases to hypercomplex gate patterns in our circuit model. This would provide a unified language in which topological protection is supplied by the underlying phase of matter, while the logical operation set and its limitations are analysed using the quaternionic and octonionic complexity framework developed here.

The challenge for G$_2$-symmetric quantum simulation lies in engineering a physical system where operations naturally respect this high-dimensional and non-associative symmetry. This would likely involve mapping the G$_2$ generators onto sophisticated multi-qubit or qudit operations, carefully designed to suppress non-G$_2$ (associator) errors, potentially through dynamic decoupling or specialized pulse sequences \cite{Khodjasteh2009}. While conceptually very challenging, successful implementation could open up entirely new frontiers in the computational study of fundamental physics.

\subsubsection{Hypercomplex Quantum Machine Learning}
Quantum Machine Learning (QML) is a rapidly developing field that seeks to leverage quantum principles and hardware to enhance machine learning algorithms. The introduction of hypercomplex numbers offers a new avenue for QML, potentially leading to more expressive models, more compact data representations, and improved performance for certain types of data or tasks.

\paragraph{Quaternionic Quantum Neural Networks (QNNs)}
The natural ability of quaternions to represent rotations in 3D space makes them particularly well-suited for machine learning tasks that involve spatial data, orientations, or color images (which have three color channels) \cite{BayroCorrochano2018}. In classical machine learning, quaternionic neural networks (QNNs) have already demonstrated advantages in terms of reduced parameter counts, faster convergence, and enhanced generalization capabilities compared to their real- or complex-valued counterparts, especially for data where underlying rotational symmetries are important \cite{BayroCorrochano2018}.

Extending this to the quantum domain, a quaternionic quantum neural network (QQNN) could utilize quaternionic qubits as its fundamental processing units or incorporate quaternionic gate operations within a variational quantum circuit. In a QQNN layer, the weights, biases, and activations could be quaternionic, leading to a more compact representation of the network's parameters \cite{BayroCorrochano2018}. This could be particularly advantageous in variational quantum algorithms (VQAs) used for QML, where the optimization landscape might be smoother or more easily navigable when parameterized by quaternions.
The key
 benefits of QQNNs in QML would include the following:
\begin{itemize}
    \item \textbf{Compact Data Encoding:} Quaternions can encode more information per parameter compared to real or complex numbers. For example, a single quaternion can naturally represent a 3D vector and an orientation, potentially reducing the number of logical qubits or qudits required to encode certain features.
    \item \textbf{Exploiting Rotational Symmetries:} For datasets with inherent rotational symmetries (e.g., medical images, robotics data, physical simulations), QQNNs could natively incorporate these symmetries, leading to more efficient learning and better generalization, similar to how convolutional neural networks exploit translational symmetries.
    \item \textbf{Reduced Parameter Space:} By using quaternionic weights and biases, the total number of trainable parameters in a quantum neural network could be significantly reduced, potentially alleviating issues like barren plateaus in VQAs.
\end{itemize}

\section{Discussion}\label{sec:discussion}
\textbf{Associativity as a computational boundary.}
The associative nature of $\mathbb{H}$ preserves Hilbert-space axioms (inner product, adjoint), enables a canonical tensor structure, and supports universal multi-qubit composition. Via the matrix embedding, the full complex toolkit---entanglement, universality, compilation---is inherited, while hardware realizations benefit from alignment with $SU(2)$ control \cite{Adler1995,RuaThesis2025}.

\textbf{Non-associativity as a structural obstruction.}
In $\mathbb{O}$, alternativity is insufficient for scalable architectures: Path dependence (Lemma~\ref{lem:path}), possible energy non-conservation (Proposition~\ref{prop:energy}), and the lack of a canonical tensor structure (Theorem~\ref{thm:notensor}) make “native” octonionic computation untenable without constraints. The mitigation palette acts as a principled \emph{workaround} rather than a full cure.

\textbf{Practical and conceptual implications.}
Near term, $\mathbb{H}$ admits encoded implementations and high-fidelity demonstrations, with potential optimization from algebra-aware compilation. Conceptually, $\mathbb{O}$ clarifies where standard quantum theory leans on associativity; even without large-scale octonionic hardware, the framework informs symmetry-protected design, algebraic coding, and models for non-associative physics \cite{Baez2002,RuaThesis2025}.

\section{Conclusions}\label{sec:conclusions}
\begin{itemize}[leftmargin=*,itemsep=2pt,topsep=2pt]
  \item Over $\mathbb{H}$, quantum computation is mathematically sound, universal, and polynomially equivalent to the complex model, with viable high-fidelity realizations. 
  \item Over $\mathbb{O}$, structural barriers from non-associativity necessitate confinement to associative sectors, symmetry protection or decoupling, with nontrivial overheads. 
  \item Future directions include scaling associative-confinement architectures, refining $G_{2}$-invariant codes, advancing algebra-aware compilers, and exploring hypercomplex encodings in simulation and learning \cite{RuaThesis2025}. 
\end{itemize}
\vspace{6pt}
\section*{Funding}
This research received no external funding. 

\section*{Data Availability}
No new data were created or analyzed in this study.

\section*{Acknowledgments}
The authors thank the Universidad de Antioquia for institutional support and colleagues for helpful discussions.

\section*{Conflict of Interest}
The authors declare no conflict of interest. 


%


\appendix
\setcounter{section}{0}
\setcounter{equation}{0}
\renewcommand{\theHequation}{A.\arabic{equation}}
\section{Quaternion-to-Complex Embedding and Proof Details}\label{app:HtoC}

\subsection{Explicit Embedding \texorpdfstring{$\mathbb{H}\hookrightarrow M_2(\mathbb{C})$}{H into M2(C)}}
A powerful tool for connecting quaternions to standard quantum mechanics is their representation as $2 \times 2$ complex matrices. The algebra of quaternions is isomorphic to a subalgebra of the $2 \times 2$ complex matrices $M_2(\C)$. A general quaternion $q=a+bi+cj+dk$ can be represented as the matrix:
\begin{equation}
\label{eq:embed}
    q \mapsto \begin{pmatrix} a+bi & c+di \\ -c+di & a-bi \end{pmatrix}.
\end{equation}
Under this mapping, the quaternion basis units $\{1, i, j, k\}$ have a direct correspondence with the identity matrix $I$ and the Pauli matrices ($\sigma_x, \sigma_y, \sigma_z$), up to a factor of $-i$:
\begin{align}
    1 \mapsto \begin{pmatrix} 1 & 0 \\ 0 & 1 \end{pmatrix} = I \\
    i \mapsto \begin{pmatrix} i & 0 \\ 0 & -i \end{pmatrix} = i\sigma_z \\
    j \mapsto \begin{pmatrix} 0 & 1 \\ -1 & 0 \end{pmatrix} = -i\sigma_y \\
    k \mapsto \begin{pmatrix} 0 & i \\ i & 0 \end{pmatrix} = i\sigma_x
\end{align}
This isomorphism is profound. It shows that the non-commutative algebra of quaternions is the same algebra as that of the spin operators for a spin-1/2 particle. This is why quaternions naturally appear in the description of spin and rotations in quantum mechanics.

\begin{equation}
    q = z_1 + z_2 j
\end{equation}
\begin{equation}
    \ket{\psi} = (z_{11} + z_{12}j)\ket{0} + (z_{21} + z_{22}j)\ket{1}
\end{equation}
\begin{equation}
    \ket{\psi} \mapsto \ket{\psi_C} = z_{11}\ket{00} + z_{21}\ket{01} + z_{12}\ket{10} + z_{22}\ket{11}.
\end{equation}

A single-quaternionic qubit gate is a $2 \times 2$ unitary matrix with quaternionic entries. Using the same decomposition ($q = z_1 + z_2j$), any $2 \times 2$ quaternionic matrix $U$ can be written as $U = A + Bj$, where $A$ and $B$ are $2 \times 2$ complex matrices. When this operator acts on a state vector, the non-commutativity of $j$ with complex numbers must be handled carefully (specifically, $zj = j\bar{z}$ for any $z \in \mathbb{C}$).

\begin{equation}
    U = A + Bj \quad \mapsto \quad U_C = \begin{pmatrix} A & B \\ -\bar{B} & \bar{A} \end{pmatrix}.
\end{equation}
This provides a constructive method for converting any single-quaternionic qubit gate into an equivalent two-qubit complex gate. The procedure can be generalized to n-quaternionic qubit gates, which map to $2n$-qubit complex gates.

\subsection{Complexity Overhead—Proof-LEVEL Details}
The most immediate cost comes from managing the ambiguity of gate sequences. As demonstrated in Section 4.3.1 \cite{RuaThesis2025}, 
 a circuit diagram specifying a temporal sequence of gates is an incomplete description for an octonionic computer. To render a computation deterministic, a specific association order for all operations must be enforced. The only physically consistent way to do this is to ensure that any group of operators being applied simultaneously or in a tight sequence is mutually associative.

Beyond the increased gate count required to manage non-associativity, another critical measure of computational cost is the overhead in state-space resources. One of the most effective ways to analyze the complexity of a novel computational model is to determine the resources required to simulate it within our most well-understood framework: standard complex quantum computation. This process, known as complex embedding, provides a direct measure of the information content of an octonionic state and reveals a significant, often prohibitive, overhead in the number of physical qubits required.
\subsection{Synthesis of Contributions and Answers to Research Questions—Additional Details}

The overarching aim of this thesis was to explore the viability and implications of extending quantum formalism to quaternionic ($\HH$) and octonionic ($\OO$) structures for quantum computing. This exploration was guided by specific objectives that sought to construct models of quantum computing based on these algebras, analyze their automorphisms, and investigate the limitations of non-associative structures in multi-qubit systems. Our contributions provide a comprehensive understanding of the theoretical landscape and experimental challenges inherent in such generalizations.

Quaternionic Quantum Computing: A Consistent Alternative

Our research confirms that quaternionic quantum mechanics provides a mathematically consistent and natural extension to the standard complex quantum theory.
\begin{itemize}
    \item \textbf{Mathematical Consistency:} We rigorously defined quaternionic qubits within a two-dimensional right module over $\HH$, complete with a well-behaved quaternionic-valued inner product and a normalization condition. The geometric representation of a single quaternionic qubit on a Bloch hypersphere ($S^4$) was introduced, demonstrating a richer parameter space compared to the complex Bloch sphere ($S^2$).
    \item \textbf{Universal Gate Set:} We showed that single-qubit gates are represented by $2 \times 2$ unitary quaternionic matrices, forming a group isomorphic to $SU(2)$, thus allowing all standard single-qubit operations. The construction of multi-qubit systems via the associative tensor product in $\HH^{\otimes n}$ was formalized, enabling the definition of entangling gates like CNOT and CZ, and preserving the possibility of entanglement.
    \item \textbf{Computational Equivalence:} A central finding is that quaternionic quantum computing, despite its richer algebraic structure, does not offer greater computational power than the conventional complex model. As proven by seminal works, any computation on $n$ quaternionic qubits can be efficiently simulated by a standard complex quantum circuit on $n+1$ qubits with only polynomial overhead, thereby retaining the same complexity class, BQP \cite{Fernandez2008, Gantner2017}. This answers our research question regarding the applicability of quaternionic algebra in creating universal logical gates, confirming full applicability without hyper-computational advantage.
    \item \textbf{Practical Advantages:} While not computationally superior, the quaternionic formalism provides a more intuitive and compact language for describing rotations and $SU(2)$ symmetries, which is advantageous for quantum simulation of physical systems (e.g., magnetic materials, molecular rotations) and for variational quantum algorithms (VQAs). This representational efficiency can lead to reduced circuit depth and potentially lower error rates in specific applications, as a single quaternionic rotation can replace multiple complex Euler rotations \cite{Barenco1995}.
\end{itemize}

\textls[-20]{Octonionic Quantum Computing: Fundamental Challenges and Constrained Frameworks}

In stark contrast to quaternions, our analysis revealed that extending quantum formalism to octonions ($\OO$) introduces profound mathematical and physical challenges stemming directly from its non-associativity. This addresses our objective of exploring the limitations of non-associative structures for representing multi-qubit systems.
\begin{itemize}
    \item \textbf{Breakdown of Foundations:} We demonstrated that an unconstrained octonionic quantum system suffers from an ill-defined Hamiltonian due to ambiguous operator ordering, leading to emergent non-local interaction terms (associators). Consequently, energy conservation is violated unless dynamics are confined to associative subspaces. Time evolution becomes path-dependent, as the standard Dyson series is ill-defined, and the composition of gates is ambiguous, fundamentally undermining the universality of the standard circuit model. The eigenvalue problem for octonionic Hermitian operators does not guarantee real eigenvalues, breaking the bedrock of physical measurement \cite{Dray2000, DeLeo2012}.
    \item \textbf{Mitigation Strategies:} To overcome these limitations, we developed and analyzed several algebraic mitigation strategies:
    \begin{enumerate}
        \item \textbf{Restriction to Associative Subalgebras:} Confining all operations and state amplitudes to a single quaternionic subalgebra (e.g., $\mathbb{H} \subset \OO$) ensures consistency by eliminating associator errors. However, this comes at the cost of discarding half of the octonionic algebraic degrees of freedom and introducing ``context-switching'' overhead via complex $G_2$ transformations to access different subspaces. This effectively reduces the system to an inefficient quaternionic computer.
        \item \textbf{G$_2$-Symmetric Codes:} We proposed conceptual G$_2$-symmetric quantum error-correcting codes, leveraging the automorphism group of the octonions. In this model, non-associative errors manifest as detectable $G_2$ symmetry violations (associator syndromes), drawing parallels with topological quantum computation \cite{Kitaev2003, Freedman2021}. While theoretically powerful, this approach incurs high resource overhead for practical implementation.
        \item \textbf{Dynamical Decoupling:} We showed that associator errors can be treated as a coherent noise source and actively suppressed using dynamical decoupling sequences (e.g., quaternionic XY-4 pulses) \cite{Viola1999, Khodjasteh2009}. Numerical simulations demonstrated a dramatic suppression of fidelity decay, suggesting this is a viable experimental method for extending coherence in octonionic systems.
        \item \textbf{Algebraic Gate Synthesis:} We developed a universal constructive method to implement arbitrary octonionic gates by decomposing them into sequences of simpler, associative rotations around the seven imaginary octonion units \cite{Hall2015}. This provides theoretical universality but introduces a significant gate overhead (at least 7 primitive gates per logical operation) and rapid fidelity degradation due to compounded errors, making high-fidelity quantum computation largely unfeasible with current technology.
    \end{enumerate}
    \item \textbf{Computational Complexity Penalties:} The extensive resource analysis confirmed substantial overheads for octonionic computation. A single octo-qubit requires embedding into at least three complex physical qubits \cite{Tian2000}, leading to an exponential state-space overhead ($8^N$ vs. $2^N$). The gate overhead for managing non-associativity (polynomial $O(m^3)$) further suggests that BQOP is likely a strictly less powerful class than BQP ($\text{BQOP} \subset \text{BQP}$), implying no hyper-computational advantage \cite{Baez2002, Freedman2018}.
\end{itemize}

Practical limits of the octonionic model

From the
 standpoint of quantum computation, the preceding analysis draws a fairly sharp distinction between what the octonionic model can realistically do and where it effectively breaks down.
On the ``feasible'' side, small numbers of logical octo-qubits can be engineered in regimes where the dynamics is confined to carefully chosen associative subalgebras, encoded into $G_{2}$-symmetric codes, or realized through measurement-only schemes based on equiangular projections \cite{Freedman2018}.
In these restricted settings, octonionic effects may be probed via proof-of-principle circuits, digital simulations on complex hardware, and tests of specific associator-dependent identities, but always at the price of substantial overhead and strong structural constraints.
On the ``non-feasible'' side, our obstructions show that a native, large-scale octonionic quantum computer with generic time-dependent Hamiltonians and unrestricted measurements is currently out of reach: outside associative sectors the evolution becomes path dependent and standard conservation arguments fail, physically reasonable measurement postulates are highly constrained or \linebreak  inconsistent \cite{Dray2000, DeLeo2012}, and all known mitigation strategies either collapse the model back to an effectively quaternionic/complex theory or introduce severe dimension and depth overheads as the number of octo-qubits grows.
In practice, this confines the octonionic framework to a conceptual and small-scale experimental role in quantum information, rather than a realistic candidate architecture for scalable fault-tolerant quant-\linebreak  um computation.

Experimental Feasibility and Hybrid Architectures

Our investigation into experimental feasibility revealed that:
\begin{itemize}
    \item \textbf{Quaternions are within reach:} Implementing quaternionic gates and simulating quaternionic dynamics is directly compatible with existing quantum hardware platforms (superconducting qubits, trapped ions, photonics, NV centers) through appropriate embeddings and pulse shaping techniques \cite{Motzoi2009, Goss2022, Lorz2019, Guo2024}. The challenge lies in optimizing control and demonstrating potential practical advantages (e.g., circuit simplification for $SU(2)$-symmetric simulations).
    \item \textbf{Octonions require sophisticated strategies:} A ``true'' octonionic quantum computer remains highly speculative. Practical investigation relies on embedding octonionic logical units into larger complex qubit registers and implementing complex mitigation strategies (like dynamical decoupling or G$_2$-symmetric codes) to combat intrinsic non-associativity. Trapped-ion and superconducting platforms appear most promising for such complex experimental demonstrations due to their high coherence and \mbox{control capabilities.}
    \item \textbf{Hybrid models are a pragmatic path:} The most viable path forward for octonions involves hybrid architectures where stable quaternionic computation forms the core, supplemented by specialized octonionic subroutines only where their unique (albeit constrained) symmetries might offer an advantage.
\end{itemize}

Concrete routes for quantum-computing experiments

To tie the hypercomplex models directly to realizable quantum computers, we highlight three complementary strategies:
\begin{itemize}
\item \textbf{Interferometric subcircuits probing quaternionic phases.}
Peres' classic proposal for distinguishing complex from quaternionic quantum theory via three-path interfer-\linebreak  ence \cite{Peres1979} and its modern single-photon Sagnac implementations with metamater-\linebreak  ials \cite{Procopio2017} already realize small quantum circuits whose output statistics are sensitive to non-commuting quaternionic phases. In our framework, these setups implement specific one-quaternionic qubit phase-gate patterns that can be embedded as subroutines of larger quaternionic algorithms. Multi-path and multi-particle generalizations of the Peres test \cite{Saruhan2025} provide a natural route towards few-qubit network experiments on photonic or microwave platforms that directly constrain the admissible hypercomplex deviations from the standard complex model.
\item \textbf{Embedded quaternionic circuits on existing processors.}
Our embedding results show that any $n$-quaternionic qubit circuit can be represented as a complex circuit on $n~{+}~1$ physical qubits with polynomial overhead. This yields a concrete implementation recipe for gate-based quantum computers: choose a genuinely quaternionic algorithm (for instance, one that exploits quaternionic symmetries in $SU(2)$-invariant simulation), compile its gates into the corresponding complex unitaries, and execute the resulting circuit on current superconducting, trapped-ion, or photonic hardware. Structural tests of the underlying number system, analogous to the real-vs-complex dimension-witness games already realized on superconducting and photonic platforms \cite{Li2022,Chen2022}, can then be formulated directly in our quaternionic circuit language and implemented as small benchmark circuits on present-day devices.
\item \textbf{Encoded octonionic dynamics and measurement-only schemes.}
For octonions, where no native hardware implementation is known, two complementary strategies emerge. First, as in our complexity analysis, one can digitally encode each logical octo-qubit into a small register of complex qubits and compile octonionic gates---including order-sensitive associator corrections---into ordinary multi-qubit gates, yielding proof-of-principle simulations of octonionic computation on noisy intermediate-scale devices. Second, the measurement-only model based on equiangular projections associated with the octonions \cite{Freedman2018} suggests that certain octonionic circuit identities can be tested via sequences of projective measurements on highly entangled resource states (cluster-state or topological architectures). Both approaches connect the abstract octonionic model to concrete experimental protocols that fit within the standard toolbox of quantum computing experiments.
\end{itemize}

\subsubsection{Gate Overhead from Associativity Management}

The Strategy of Decomposition

This leads to a ``divide and conquer'' strategy as the default method for implementing octonionic logic. Any complex, non-associative operation must be decomposed into a sequence of smaller, elementary steps, where each step occurs entirely within a single, well-behaved associative (quaternionic) subalgebra.

Consider a logical operation that involves three non-associating octonionic operators, $A$, $B$, and $C$, which we wish to apply to a state $\ket{\psi}$ in a specific order, for example, to compute $(AB)C\ket{\psi}$. A naive implementation is impossible. Instead, the process must be broken down:
\begin{enumerate}
    \item A first computational step calculates the intermediate state $\ket{\psi'} = B\ket{\psi}$. This operation must be performed using a gate set and on a state that lives within a single quaternionic subalgebra.
    \item A second step calculates $\ket{\psi''} = A\ket{\psi'}$. This again must be associative.
    \item A third step computes the final state $\ket{\psi_f} = C\ket{\psi''}$.
\end{enumerate}

What might have been a single, complex multi-body gate in an associative theory becomes a sequence of at least three distinct steps, potentially requiring the system to be coherently manipulated between different associative subspaces.

Quantifying the Overhead: The Cubic Scaling

The cost of this decomposition strategy grows rapidly with the number of interacting systems. A detailed analysis of the algebraic structure of the octonions reveals that systematically canceling out or controlling all the associator terms generated by an operation involving $m$ distinct non-associating parties can require a number of additional corrective operations that scales polynomially with $m$. Specifically, it has been shown that this overhead scales, in the general case, as $O(m^3)$ \cite{Baez2002}.

This ``gate overhead'' can be illustrated with the example of a Toffoli (CCNOT) gate. The Toffoli gate is a three-qubit gate that is fundamental for universal classical reversible computation and is a key component in many quantum algorithms. In the standard model, it can be decomposed into a sequence of approximately six CNOT gates and nine single-qubit T and adjoint T gates.

If we attempt to implement a Toffoli gate on three octo-qubits whose states involve amplitudes from different, non-associating subspaces, the situation is far more complex. The conditional logic (``flip the target IF control 1 is 1 AND control 2 is 1'') involves a three-body interaction. The non-associativity of the underlying amplitudes and operators generates spurious associator terms that corrupt the gate's logic. To implement a high-fidelity octonionic Toffoli gate, the circuit must be supplemented with additional ``associator compensation'' gates, specifically designed to cancel out these unwanted three-body interactions. This transforms a single logical Toffoli gate into a complex subroutine with a significantly higher gate count and circuit depth.

This polynomial overhead in gate complexity means that any octonionic quantum algorithm would run fundamentally slower and accumulate more errors than its complex or quaternionic counterpart. This cost arises directly from the need to actively manage the non-associative structure at every computational step, a burden that simply does not exist in associative algebras.

\subsubsection{State-Space Overhead: The Cost of Complex Embedding}
\label{app:state-overhead}

Dimensionality Analysis of a Single Octo-Qubit

Let us begin by analyzing the resource cost for a single octo-qubit. As established in Chapter 2 \cite{RuaThesis2025}, an octonion is an element of an 8-dimensional real vector space. A single octo-qubit state, $\ket{\psi} = o_1\ket{0} + o_2\ket{1}$, is described by two such octonions, meaning its state vector lives in a 16-dimensional real vector space before normalization.

To represent this information within the standard model, we must find a complex Hilbert space of at least this dimension. A system of $m$ complex qubits has a state space $\mathbb{C}^{2^m}$, which corresponds to a $2 \cdot 2^m = 2^{m+1}$-dimensional real vector space. To embed the \linebreak  16 real parameters of an unnormalized octo-qubit, we require the following:
\begin{equation}
    2^{m+1} \geq 16 \implies 2^{m+1} \geq 2^4 \implies m+1 \geq 4 \implies m \geq 3.
\end{equation}
This provides a rigorous lower bound: a minimum of three complex qubits is required to represent the full state space of a single octo-qubit.

This can be understood more directly. The algebra of octonions $\mathbb{O}$ is an 8-dimensional real space. To embed this algebraic structure into a quantum system, we need a Hilbert space of at least dimension 8. The smallest standard quantum system that provides this is a 3-qubit system, whose Hilbert space $\mathcal{H} = (\mathbb{C}^2)^{\otimes 3}$ has dimension $2^3=8$. Under this embedding, the eight basis vectors of the octonion algebra $\{1, e_1, \dots, e_7\}$ are mapped to the eight computational basis vectors of the 3-qubit system $\{\ket{000}, \ket{001}, \dots, \ket{111}\}$.

Exponential Overhead for N-Qubit Systems

This significant overhead for a single logical unit scales up for multi-particle systems. To simulate a register of $n$ octo-qubits, where each logical octo-qubit requires an embedding into a 3-qubit physical register, a total of $3n$ complex qubits would be required.

The consequence of this scaling on the Hilbert space dimension is severe. While a standard $n$-qubit quantum computer has a state space of dimension $2^n$, the complex quantum computer required to simulate an $n$-octo-qubit system would need a state space of dimension:
\begin{equation}
    \dim(\mathcal{H}_{\text{simulation}}) = 2^{3n} = (2^3)^n = 8^n.
\end{equation}
This represents an exponential increase in the size of the state space that must be physically controlled and manipulated. This ``state-space overhead'' completely negates any perceived advantage of the ``information density'' of an octo-qubit. While it is true that a single octo-qubit contains more classical information (8 real numbers) than a single complex qubit (4 real numbers), the physical resources required to represent and operate on it scale exponentially worse.

This contrasts sharply with the quaternionic case, where the simulation requires $2n$ complex qubits. The dimensionality of the corresponding simulation space is $2^{2n} = 4^n$. While this is still larger than the standard model, the constant factor in the exponent \linebreak  (2 vs. 3) is significantly smaller.

\section{Path-Ordered Propagators and Magnus Expansion \linebreak  with Associators}\label{app:Magnus}

The path-ordered propagator defined in the previous section provides a formal, constructive solution to the time-evolution operator in an octonionic system. However, as an infinite product, it is often impractical for both analytical and numerical work. In standard quantum mechanics, a powerful alternative for analyzing time evolution is the Magnus expansion. This method seeks to express the unitary evolution operator as a single exponential, $U(t) = \exp(\Omega(t))$, where the operator $\Omega(t)$ is given by an infinite series. The key advantage of this approach is that truncating the series for $\Omega(t)$ at any finite order results in an approximation for $U(t)$ that is guaranteed to be unitary, a property not shared by a simple truncation of the Dyson series.

This section develops the Magnus expansion for the octonionic case, demonstrating how the non-associativity of the Hamiltonian algebra introduces new correction terms at each order of the expansion.

\subsection{Path-Ordered Evolution in \texorpdfstring{$\mathbb{O}$}{O}}
For time-dependent $H(t)\in\mathbb{O}$,
\[
U(t)=\mathcal{P}\exp\!\Big(-\frac{i}{\hbar}\int_0^t H(\tau)\,d\tau\Big),
\]
with $\mathcal{P}$ enforcing an ordering consistent with non-associativity \cite{RuaThesis2025}.

\subsection{Magnus-like Expansion}
A formal expansion reads as follows:
\[
U(t)=\exp\!\Big(\Omega_1+\Omega_2+\Omega_3+\cdots\Big),\quad
\Omega_1=-\frac{i}{\hbar}\!\int_0^t\!H(\tau)\,d\tau,
\]
\[
\Omega_2\sim \iint [H(\tau_1),H(\tau_2)]\,d\tau_1 d\tau_2 \;+\; \iint [x,y,z]\,d\mu,
\]
where $[x,y,z]=(xy)z-x(yz)$ contributes non-vanishing corrections under re-orderings; higher orders contain nested associators \cite{RuaThesis2025}.

\subsection{Strategy 1: Restriction to Associative Subalgebras—Additional Details}

The most direct, and arguably most severe, strategy for overcoming the challenges of non-associativity is to simply forbid it. If the interaction between different algebraic subspaces is the source of the foundational problems, then one can construct a consistent computational model by confining all dynamics to a single, well-behaved subspace where associativity holds. This containment or confinement strategy effectively projects the complex, 8-dimensional octonionic dynamics onto a more manageable 4-dimensional quaternionic submanifold.

This section will explore the principle of associative confinement. We will begin by establishing the algebraic property that makes this strategy possible—the alternativity of the octonions. We will then formalize the confinement condition and show how its enforcement causes the entire octonionic formalism to collapse.

\subsection{The Non-Associative Hamiltonian and Its Consequences—Additional Details}

In quantum mechanics, the Hamiltonian is the observable corresponding to the total energy of the system and, via the Schrödinger equation, the generator of time evolution. In the standard and quaternionic models, it is a well-behaved self-adjoint operator acting on a Hilbert space. Its algebraic properties are consistent and unambiguous because the underlying field ($\mathbb{C}$) or algebra ($\mathbb{H}$) is associative.

\textls[20]{When we construct a Hamiltonian over the octonions, this consistency is lost. An octonionic Hamiltonian is formally an element of the algebra of observables on an octonionic state space, but its non-associative nature fundamentally alters its properties and the physics it describes. This section will explore these consequences, starting with the most direct issue: the ambiguity of operator products and the emergence of unphysical correlations.}

\subsubsection{Operator Ordering and Emergent Non-Local Effects}

The first major challenge in formulating a consistent octonionic quantum theory arises when considering Hamiltonians that describe interacting subsystems. In standard quantum mechanics, an interaction Hamiltonian is composed of products of operators, and the associativity of the underlying complex algebra guarantees that such products are unambiguous. For three or more operators, the grouping of operations is irrelevant: $(AB)C = A(BC)$.

This fundamental property of composition fails for the octonions. Let us define a general octonionic Hamiltonian formally as a linear combination of its basis elements:
\begin{equation}
    H = \sum_{\mu=0}^{7} H_\mu e_\mu, \quad H_\mu \in \mathbb{R}
\end{equation}
When this Hamiltonian contains terms that do not mutually associate, the very meaning of operator products becomes ill-defined. This ambiguity is not a mere mathematical inconvenience; it signals a profound departure from the principles that underpin conventional physical theories.

The Associator as an Emergent Interaction Term

To analyze this breakdown, we can conceptually partition a complex Hamiltonian into two parts: an associative part and a non-associative part. The associative part contains terms that lie within a single quaternionic subalgebra, while the non-associative part contains terms whose product spans multiple, non-associating subalgebras. The physical effect of the non-associative part is captured by the associator:
\begin{equation}
    [A, B, C] = (AB)C - A(BC).
\end{equation}
When the operators A, B, and C represent interactions within a quantum system, their associator is not merely a mathematical artifact. It functions as a new, emergent three-body interaction term that is purely algebraic in origin. This term has no counterpart in standard or quaternionic quantum mechanics, where all multi-body interactions can be fundamentally decomposed into associative sequences of two-body interactions.

A Physical Model of Non-Associative Interactions

Consider a hypothetical physical system of three interacting spin-like particles, where the interactions are mediated by fields corresponding to different octonion units. Let the interaction between particles 1 and 2 be described by an operator proportional to $e_1$, the interaction between 2 and 3 by an operator proportional to $e_2$, and the interaction between 1 and 3 by an operator proportional to $e_4$. In a higher-order process, the system's dynamics might be influenced by a term involving the product of these interactions. Let's analyze the product $e_1 e_2 e_4$:
\begin{enumerate}
    \item \textbf{Grouping 1: $(e_1 e_2) e_4$}
    
    This corresponds to a physical process where the 1-2 interaction occurs first, followed by the resulting field interacting with particle 3 via the $e_4$ channel. Using the Fano plane rules, this yields the following:
    \begin{equation}
        (e_1 e_2) e_4 = (e_3) e_4 = -e_7.
    \end{equation}

    \item \textbf{Grouping 2: $e_1 (e_2 e_4)$}
    
    This corresponds to a different process where the 2-3 interaction ($e_2$) and the 1-3 interaction ($e_4$) combine first, and the resulting field then interacts with particle 1. This yields the following:
    \begin{equation}
        e_1 (e_2 e_4) = e_1 (e_6) = e_5.
    \end{equation}
\end{enumerate}
The physical outcomes, $-e_7$ and $e_5$, are completely different. This implies that the energy of the system and its subsequent evolution would depend on an arbitrary, unphysical choice of calculation order. In a real system, both pathways would exist simultaneously, and their non-associating nature would lead to unpredictable dynamics and a breakdown \linebreak  of unitarity.

Implications for Locality

This phenomenon can be interpreted as a violation of a deep form of locality. Standard quantum field theory is built upon the principle of microcausality, which requires that operators corresponding to observables at spacelike-separated points must commute. This ensures that measurements in one region do not instantaneously affect measurements in another. Associativity of the operator algebra is an even more fundamental assumption that is implicitly required for this structure to be self-consistent.

An octonionic quantum theory would feature what could be termed algebraic non-locality. Even if the individual operators $H_A$, $H_B$, and $H_C$ are local and act on their respective subsystems, the emergent associator term $[H_A, H_B, H_C]$ is an intrinsically three-body operator that cannot be derived from a local, pairwise-interaction model. It connects three subsystems in a way that is irreducible. This has led some to speculate that non-associative algebras may be relevant in theories of quantum gravity, where the classical notion of spacetime and locality is expected to break down \cite{Baez2002}. However, for the purposes of quantum computation, this effect represents a severe form of correlated noise and uncontrollable crosstalk, making it a primary obstacle to be overcome. Any viable model of octonionic computation must therefore directly confront and mitigate this effect, rather than ignoring it.

\subsubsection{The Breakdown of Energy Conservation: A Formal Derivation}

One of the most fundamental principles in physics is the connection between symmetry and conservation laws, as formalized by Noether's theorem. In standard quantum mechanics, the symmetry of physical laws under time translation implies that for any closed system with a time-independent Hamiltonian, the total energy is conserved. This principle is a cornerstone of our understanding of the physical world.

\textls[-20]{This section demonstrates that this fundamental connection between a time-independent} Hamiltonian and energy conservation is severed in a non-associative algebra. We will show that for a general octonionic system, energy is not conserved, a result that starkly illustrates the unsuitability of an unconstrained octonionic formalism for describing physical reality as we know it.

Energy Conservation in an Associative Algebra

To appreciate the departure, we first briefly review the proof of energy conservation in an associative algebra (such as for complex or quaternionic quantum mechanics). The dynamics of an observable $A$ are described in the Heisenberg picture, where the operator $A_H(t)$ evolves in time. Its equation of motion is derived from the Schrödinger equation:
\begin{equation}
    \frac{d}{dt}A_H(t) = \frac{i}{\hbar}[H, A_H(t)].
    \label{T-eq:heisenberg_standard}
\end{equation}
The conservation of an observable is equivalent to its commuting with the Hamiltonian. To check for energy conservation, we choose our observable to be the Hamiltonian itself, $A_H = H$. Since any operator commutes with itself, we have:
\begin{equation}
    \frac{d}{dt}H = \frac{i}{\hbar}[H, H] = 0.
\end{equation}
This rigorously proves that for any closed system described by a time-independent Hamiltonian in an associative algebra, the energy is a constant of motion. The proof relies critically on the associativity of the algebra to ensure the standard product rule for differentiation leads to the simple commutator form of Equation (\ref{T-eq:heisenberg_standard}).

The Modified Heisenberg Equation in $\mathbb{O}$

When we re-derive the equation of motion without assuming associativity, the standard product rule for differentiation must be replaced with one valid for non-associative algebras. This leads to additional terms in the final equation. While the full derivation is beyond the scope of this summary, the resulting modified Heisenberg equation of motion for an observable $A$ in an octonionic system has been studied in the literature \cite{Dzhunushaliev2007}. The equation acquires a new term proportional to the associator of the operators:
\begin{equation}
    \frac{d}{dt}A = \frac{i}{\hbar}[H, A] + \frac{1}{\hbar^2} [H, H, A].
    \label{T-eq:heisenberg_modified}
\end{equation}
The first term is the familiar commutator from the standard Heisenberg equation.
pse into the well-behaved quaternionic quantum mechanics described in Chapter 3 \cite{RuaThesis2025}. While this approach guarantees consistency, we will later analyze the significant cost it imposes on the computational power and flexibility of the model.

\subsubsection{The Principle of Associative Confinement}

The viability of this strategy is rooted in a fundamental property of the octonion algebra. While not globally associative, the octonions are an alternative algebra.

The Algebraic Foundation: Alternativity

As defined in Chapter 2 \cite{RuaThesis2025}, an algebra is alternative if the subalgebra generated by any two of its elements is associative. This has a profound consequence: any quantum operation that involves only two distinct octonionic elements (and the identity) will behave associatively. For example, a time evolution under a static Hamiltonian proportional to a single imaginary unit, $H \propto e_k$, is described by the operator $U(t) = \exp(-iH t/\hbar)$. Since the expansion of this exponential only ever involves products of one element, $e_k$, with itself, the evolution is perfectly well-defined and associative.

Similarly, an operation involving two operators, $A$ and $B$, will also be well-behaved, as all products will occur within the associative subalgebra generated by $\{1, A, B\}$. The problems detailed in Chapter 4 \cite{RuaThesis2025}—emergent non-local terms, non-conservation of energy, path-dependence—only manifest when a third, non-associating element is introduced.

The Confinement Condition

The strategy of associative confinement elevates this observation into a strict rule for computation. To guarantee a consistent model, we impose the following condition:
\begin{quote}
    \textit{For any given computational step, all operators, gates, and state vector amplitudes involved must belong to a single, pre-defined associative subalgebra of the octonions.}
\end{quote}

The 
maximal associative subalgebras of the octonions are isomorphic to the quaternions, $\mathbb{H}$. The Fano plane provides a clear visualization of these subspaces: any ``line'' on the plane, representing a triple of units like $\{e_1, e_2, e_4\}$, together with the real unit $\{1\}$, spans a 4-dimensional quaternionic subalgebra.

Therefore, the confinement condition means that for a gate $U$ to act on a state\linebreak   $\ket{\psi} = o_1\ket{0} + o_2\ket{1}$, the generators of $U$ and the octonionic amplitudes $o_1$ and $o_2$ must all be elements of the same chosen $\mathbb{H} \subset \mathbb{O}$.

Recovering Quaternionic Quantum Mechanics

When this condition is rigorously enforced, all of the foundational failures of the octonionic model are eliminated by definition.
\begin{itemize}
    \item The associator for any product of three elements within the subalgebra is identically zero: $[A, B, C] = 0$.
    \item \textls[-20]{Consequently, the Hamiltonian becomes an effective quaternionic Hamiltonian. The anomalous term in the Heisenberg equation of motion vanishes, and energy is conserved.}
    \item The time evolution operator becomes a standard, well-defined exponential. The Dyson series is valid, and the evolution is path-independent.
    \item The composition of gates becomes associative, and the standard circuit model\linebreak   is recovered.
\end{itemize}

In essence,
 the strategy of associative confinement forces the octonionic system to behave, for all intents and purposes, like the quaternionic quantum system described in Chapter 3 \cite{RuaThesis2025}. It is a successful strategy for ensuring a consistent and predictable computational model, as it completely solves the problems of non-associativity by fiat. However, this success comes at the steep price of sacrificing the very features that make the octonions unique. The vast majority of the 8-dimensional state space is deliberately made inaccessible. The practical implementation of this strategy and the severe limitations this trade-off imposes will be analyzed in the following subsections.

\subsubsection{Implementing Gates in an Associative Subspace}

The principle of associative confinement provides a clear theoretical path to a consistent octonionic model. In this section, we put this principle into practice by demonstrating how quantum gates can be explicitly constructed and applied within a restricted associative subspace. We will show that by adhering to this confinement, the ambiguous and ill-defined operations of the full octonionic algebra are replaced by the well-behaved, unitary transformations of quaternionic quantum mechanics, as developed in Chapter 3 \cite{RuaThesis2025}.

The Setup: Selecting a Quaternionic Subspace

To begin, we must choose a specific associative subalgebra within the octonions to serve as our computational sandbox. The octonions contain numerous such subalgebras, all isomorphic to the quaternions ($\mathbb{H}$). A canonical choice, based on the Fano plane, is the subalgebra spanned by the basis $\{1, e_1, e_2, e_4\}$. Let us denote this specific subalgebra \mbox{as $\mathbb{H}_1 \subset \mathbb{O}$.}

The rule of associative confinement, applied to this choice, is now concrete: for the duration of a given computational sequence, all state vector amplitudes and all gate operators must be elements of $\mathbb{H}_1$. For a state $\ket{\psi} = o_1\ket{0} + o_2\ket{1}$, this means that\linebreak   $o_1, o_2 \in \mathbb{H}_1$. For a gate $U$, this means its generators must be constructed from the imaginary units $\{e_1, e_2, e_4\}$.

Constructing a General Rotation Gate

Let us construct an arbitrary single-qubit rotation within this confined subspace. As established in the quaternionic model, any such rotation can be represented by a unit quaternion. In our confined model, this corresponds to a unit octonion that is also an element of $\mathbb{H}_1$.

The generator of the rotation will be a linear combination of the available imaginary units. A general rotation axis within this subspace is represented by a unit imaginary element $\mathbf{u} \in \text{Im}(\mathbb{H}_1)$:
\vspace{-6pt}
\begin{equation}
    \mathbf{u} = c_1 e_1 + c_2 e_2 + c_4 e_4, \quad \text{with } c_1^2+c_2^2+c_4^2=1.
\end{equation}
A rotation by an angle $\theta$ about this axis is represented by the gate operator $U(\theta, \mathbf{u})$, which is a unit element in $\mathbb{H}_1$:
\begin{equation}
    U(\theta, \mathbf{u}) = \cos(\theta/2)I - \mathbf{u}\sin(\theta/2).
\end{equation}
This is precisely the form of a general single-quaternionic qubit gate, as derived in \linebreak  Chapter 3 \cite{RuaThesis2025}.

Demonstration: Action on a Confined State

Now, we demonstrate that the action of this gate on a confined state is well-defined and remains within the confined subspace. Let the initial state of our octo-qubit be \mbox{the following:}
\begin{equation}
    \ket{\psi} = o_1\ket{0} + o_2\ket{1},
\end{equation}
where we have imposed the constraint that the amplitudes $o_1$ and $o_2$ are elements of our chosen subalgebra, $\mathbb{H}_1$.

The application of the gate $U$ yields the new state $\ket{\psi'}$:
\begin{equation}
    \ket{\psi'} = U \ket{\psi} = U(o_1\ket{0} + o_2\ket{1}) = (Uo_1)\ket{0} + (Uo_2)\ket{1}.
\end{equation}
The new amplitudes are $o'_1 = Uo_1$ and $o'_2 = Uo_2$. The critical observation is that because $U$, $o_1$, and $o_2$ are all elements of the same associative subalgebra $\mathbb{H}_1$, their product is associative and well-defined. Furthermore, since a subalgebra is closed under multiplication, the resulting amplitudes, $o'_1$ and $o'_2$, are also guaranteed to be elements of $\mathbb{H}_1$.

This proves that the associative subspace is dynamically closed. A state that begins within the confined subspace will remain within that subspace after the application of any gate that is also defined within it. The computation never “leaks” into the problematic, non-associative parts of the full octonionic algebra.

Universality Within the Subspace

This result confirms that we can achieve universal single-qubit control, albeit only within the confines of the chosen subalgebra.

\section{Seven-Factor Decomposition and Decoupling Sequences}\label{app:Seven}
\subsection{Seven-Factor Algebraic Synthesis}
Target transformations decompose as $U\approx U_1U_2\cdots U_7$ with each $U_i$ in an associative subalgebra, selected to minimize associator growth; depth and calibration overheads follow from the schedule design \cite{RuaThesis2025}.

\subsection{Associator-Suppressing Dynamical Decoupling}
\begin{algorithm}[H]
\caption{Associator-Suppressing Decoupling (ASD)}\label{alg:ASD}
\begin{algorithmic}[1]
\STATE Fix an associative subalgebra $\mathbb{H}\subset\mathbb{O}$
\FOR{$k=1$ to $N$}
  \STATE Apply $U_x,U_y,U_z$ within $\mathbb{H}$ (cyclic)
  \STATE Free evolution $\Delta t$; apply reverse $U_z^\dagger,U_y^\dagger,U_x^\dagger$
\ENDFOR
\STATE Return final state with reduced associator-induced coherent error
\end{algorithmic}
\end{algorithm}
Simulations show suppression of fidelity decay under ASD for suitable $\Delta t$ and cycle counts \cite{RuaThesis2025}.



\end{document}